\newif\ifpagetitre            \pagetitretrue
\newtoks\hautpagetitre        \hautpagetitre={\hfil}
\newtoks\baspagetitre         \baspagetitre={\hfil}
\newtoks\auteurcourant        \auteurcourant={\hfil}
\newtoks\titrecourant         \titrecourant={\hfil}

\newtoks\hautpagegauche       \newtoks\hautpagedroite
\hautpagegauche={\hfil\the\auteurcourant\hfil}
\hautpagedroite={\hfil\the\titrecourant\hfil}

\newtoks\baspagegauche \baspagegauche={\hfil\tenrm\folio\hfil}
\newtoks\baspagedroite \baspagedroite={\hfil\tenrm\folio\hfil}

\headline={\ifpagetitre\the\hautpagetitre
\else\ifodd\pageno\the\hautpagedroite
\else\the\hautpagegauche\fi\fi}

\footline={\ifpagetitre\the\baspagetitre
\global\pagetitrefalse
\else\ifodd\pageno\the\baspagedroite
\else\the\baspagegauche\fi\fi}

\vsize=9.0in\voffset=1cm
\looseness=2


\message{fonts,}

\font\tenrm=cmr10
\font\ninerm=cmr9
\font\eightrm=cmr8
\font\teni=cmmi10
\font\ninei=cmmi9
\font\eighti=cmmi8
\font\ninesy=cmsy9
\font\tensy=cmsy10
\font\eightsy=cmsy8
\font\tenbf=cmbx10
\font\ninebf=cmbx9
\font\tentt=cmtt10
\font\ninett=cmtt9

\font\ninesl=cmsl9
\font\eightsl=cmsl8

\font\nineit=cmti9
\font\eightit=cmti8

\skewchar\ninei='177 \skewchar\eighti='177
\skewchar\ninesy='60 \skewchar\eightsy='60

\def\eightpoint{\def\rm{\fam0\eightrm} 
\normalbaselineskip=9pt
\normallineskiplimit=-1pt
\normallineskip=0pt

\textfont0=\eightrm \scriptfont0=\sevenrm \scriptscriptfont0=\fiverm
\textfont1=\ninei \scriptfont1=\seveni \scriptscriptfont1=\fivei
\textfont2=\ninesy \scriptfont2=\sevensy \scriptscriptfont2=\fivesy
\textfont3=\tenex \scriptfont3=\tenex \scriptscriptfont3=\tenex
\textfont\itfam=\eightit  \def\it{\fam\itfam\eightit} 
\textfont\slfam=\eightsl \def\sl{\fam\slfam\eightsl} 

\setbox\strutbox=\hbox{\vrule height6pt depth2pt width0pt}%
\normalbaselines \rm}

\def\ninepoint{\def\rm{\fam0\ninerm} 
\textfont0=\ninerm \scriptfont0=\sevenrm \scriptscriptfont0=\fiverm
\textfont1=\ninei \scriptfont1=\seveni \scriptscriptfont1=\fivei
\textfont2=\ninesy \scriptfont2=\sevensy \scriptscriptfont2=\fivesy
\textfont3=\tenex \scriptfont3=\tenex \scriptscriptfont3=\tenex
\textfont\itfam=\nineit  \def\it{\fam\itfam\nineit} 
\textfont\slfam=\ninesl \def\sl{\fam\slfam\ninesl} 
\textfont\bffam=\ninebf \scriptfont\bffam=\sevenbf
\scriptscriptfont\bffam=\fivebf \def\bf{\fam\bffam\ninebf} 
\textfont\ttfam=\ninett \def\tt{\fam\ttfam\ninett} 

\normalbaselineskip=11pt
\setbox\strutbox=\hbox{\vrule height8pt depth3pt width0pt}%
\let \smc=\sevenrm \let\big=\ninebig \normalbaselines
\parindent=1em
\rm}

\def\tenpoint{\def\rm{\fam0\tenrm} 
\textfont0=\tenrm \scriptfont0=\ninerm \scriptscriptfont0=\fiverm
\textfont1=\teni \scriptfont1=\seveni \scriptscriptfont1=\fivei
\textfont2=\tensy \scriptfont2=\sevensy \scriptscriptfont2=\fivesy
\textfont3=\tenex \scriptfont3=\tenex \scriptscriptfont3=\tenex
\textfont\itfam=\nineit  \def\it{\fam\itfam\nineit} 
\textfont\slfam=\ninesl \def\sl{\fam\slfam\ninesl} 
\textfont\bffam=\ninebf \scriptfont\bffam=\sevenbf
\scriptscriptfont\bffam=\fivebf \def\bf{\fam\bffam\tenbf} 
\textfont\ttfam=\tentt \def\tt{\fam\ttfam\tentt} 

\normalbaselineskip=11pt
\setbox\strutbox=\hbox{\vrule height8pt depth3pt width0pt}%
\let \smc=\sevenrm \let\big=\ninebig \normalbaselines
\parindent=1em
\rm}

\message{fin format jgr}

\hautpagegauche={\hfill\ninerm\the\auteurcourant}
\hautpagedroite={\ninerm\the\titrecourant\hfill}
\auteurcourant={R.G.Novikov}
\titrecourant={An effectivization of the global reconstruction in the
Gel'fand-Calderon inverse problem in 3D}

\magnification=1200
\font\Bbb=msbm10
\def\C{\hbox{\Bbb C}}
\def\R{\hbox{\Bbb R}}
\def\S{\hbox{\Bbb S}}
\def\N{\hbox{\Bbb N}}
\def\b{\backslash}
\def\v{\varphi}
\def\pa{\partial}
\def\ep{\varepsilon}

\vskip 2 mm
\centerline{\bf An effectivization of the  global reconstruction}
\centerline{\bf in the
Gel'fand-Calderon inverse problem in three dimensions}
\vskip 2 mm
\centerline{\bf R. G. Novikov}
\vskip 2 mm
\centerline{18 October 2008}
\vskip 2 mm

\noindent
{\ninerm CNRS, Laboratoire de Math\'ematiques Jean Leray (UMR 6629),
Universit\'e de Nantes, BP 92208,}

\noindent
{\ninerm F-44322, Nantes cedex 03, France}

\noindent
{\ninerm e-mail: novikov@math.univ-nantes.fr}
\vskip 4 mm
{\bf Abstract.}
By developing the $\bar\pa$-approach to global "inverse scattering" at zero
energy we give a principal effectivization of the global reconstruction
method for the
Gel'fand-Calderon inverse boundary value problem in three dimensions. This
work goes back to results published by the author in 1987, 1988 and proceeds
from recent progress in the $\bar\pa$-approach to inverse scattering in 3D
published by the author in 2005, 2006.

\vskip 2 mm
{\bf Mathematics Subject Classification (2000):} 35R30 (primary); 81U40, 
86A20 (secondary).

\vskip 2 mm
{\bf 1. Introduction}

The global reconstruction method for the
Gel'fand-Calderon inverse boundary value problem in 3D was proposed  for the
first time  by the author in 1987, 1988, see [HN] (Note added in proof) and
[No1]. The sheme of this global reconstruction  can be presented as follows:
$$\Phi \to h\big|_{\bar\Theta_{\rho}} \to
v_{\rho} \to v,\eqno(*)$$
where $\Phi$ denotes boundary measurements, $h$ is the Faddeev
generalized scattering amplitude (at zero energy for simplicity) defined
in complex domain $\Theta$, $\bar\Theta_{\rho}$ is the subset of $\Theta$ with
 the imaginary part not greater than $\rho$, $v_{\rho}$ is an approximate
inverse scattering reconstruction from $h\big|_{\bar\Theta_{\rho}}$,
$v_{\rho}\to v$ for $\rho\to +\infty$, where $v$ is the unknown potential
to be reconstructed from $\Phi$. The global uniqueness in the Gel'fand-
Calderon problem in 3D follows from this global reconstruction as a
corollary, see [No1].

Slightly earlier with respect to [No1] the global uniqueness in the
Calderon problem in 3D was proved in [SU]. Note that [SU] gives no
reconstruction method.

Slightly later with respect to [No1] the global reconstruction from boundary
 measurements in 3D similar to the global reconstruction of [No1] was also
published in [Na1]. Note that [Na1] contains a reference to the preprint of
[No1].

The Gel'fand-Calderon problem is formulated as Problem 1.2 of Subsection 1.2
of this introduction.

Reconstruction problems from the Faddeev generalized scattering amplitude $h$
in the complex domain at zero energy are formulated as Problems 1.1a, 1.1b
and 1.1c of Subsection 1.1 of this introduction.

The reduction of Problem 1.2 to Problems 1.1 is given by formulas and
equations (1.23)-(1.25) mentioned in Subsection 1.2 of this introduction,
see also [No2] for more advanced version of these formulas and equations.

For a long time the main disadvantage of the global reconstruction (*) in 3D
was related with the  following two facts:

(1) The determination of $h\big|_{\bar\Theta_{\rho}}$ from $\Phi$ via
formulas and equations of the type (1.23)-(1.25) is stable for relatively
small $\rho$, but is very unstable for $\rho\to +\infty$ in the points of
$\bar\Theta_{\rho}$ with sufficiently great imaginary part, see Subsection
1.2 of this introduction.

(2) The decay of the error $v-v_{\rho}$ for $\rho\to +\infty$ was very slow
(not faster than $O(\rho^{-1})$ even for infinitely smooth compactly supported
 $v$) in existing global  results for stable construction of $v_{\rho}$
from $h\big|_{\bar\Theta_{\rho}}$ in 3D, see Remarks 1.1, 1.2, 1.3 of
Subsection 1.1 of this introduction.

As a corollary, the global reconstruction (*) in 3D was not efficient with
respect to its stability properties.

The key point is that in the present work we give a global and stable
construction of $v_{\rho}$ from  $h\big|_{b\Theta_{\rho}}$ in 3D, where
$b\Theta_{\rho}$ denotes the boundary of $\bar\Theta_{\rho}$, with rapid
decay of the error $v-v_{\rho}$ for $\rho\to +\infty$ (in particular, with
$v-v_{\rho}=O(\rho^{-\infty})$ for $v$ of the Schwartz class). This gives
a principal effectivization of the global reconstruction (*) with respect to
its stability properties.

Our new results are presented in detail below in Subsections 1.1, 1.2, 1.3
(of the introduction) and in Sections 2 and 6. These results were obtained
proceeding from [No3], [No4].

\vskip 2 mm
{\it 1.1. Inverse scattering at zero energy.}
Consider the equation
$$-\Delta\psi+v(x)\psi=0,\ \ x\in\R^d,\ \ d\ge 2,\eqno(1.1)$$
where
$$\eqalign{
&v\ \ {\rm is\ a\ sufficiently\ regular\ function\ on}\ \ \R^d \cr
&{\rm with\ sufficient\ decay\ at\ infinity} \cr}\eqno(1.2)$$
(precise assumptions on $v$ are specified below in this introduction and
 in Section 2).

Equation (1.1) arises, in particular, in quantum mechanics, acoustics,
electrodynamics. Formally, (1.1) looks as the Schr\"odinger equation with
potential $v$ at fixed energy $E=0$.

For equation (1.1), under assumptions (1.2), we consider the
Faddeev generalized scattering amplitude $h(k,l)$, where $(k,l)\in\Theta$,
$$\Theta=\{k\in\C^d,\ l\in\C^d:\ k^2=l^2=0,\ \ Im\,k=Im\,l\}.\eqno(1.3)$$
Given $v$, to determine $h$ on $\Theta$ one can use, in
particular, the formula
$$h(k,l)=H(k,k-l),\ \ (k,l)\in\Theta,\eqno(1.4)$$
and the linear integral equation
$$H(k,p)=\hat v(p)-\int\limits_{\R^d}
{\hat v(p+\xi)H(k,-\xi)d\xi\over {\xi^2+2k\xi}},\ \ k\in\Sigma,\ \
p\in\R^d,\eqno(1.5)$$
where
$$\eqalignno{
&\hat v(p)=(2\pi)^{-d}\int\limits_{\R^d}e^{ipx}v(x)dx,\ \ p\in\R^d,&(1.6)\cr
&\Sigma=\{k\in\C^d:\ k^2=0\}.&(1.7)\cr}$$

For more details concerning definitions of $h$, see [HN, Section 2.2],
[No1, Section 2] and [No4, Sections 1 and 3].

Actually, $h$ on $\Theta$ is a zero energy restriction of a function $h$
introduced by Faddeev as an extension to the complex domain of the
classical scattering amplitude for the Schr\"odinger equation at positive
energies (see [F2], [HN]).
Note that the restriction $h\big|_{\Theta}$ was not considered in
Faddeev's works and that $h$ in its zero energy
restriction was considered for the first time in [BC] for $d=3$ in the
framework of Problem 1.1a formulated below. The Faddeev function $h$ was,
actually, rediscovered in [BC]. The fact that $\bar\pa$-scattering data of
[BC] coincide with the Faddeev function $h$ was observed, in particular,
in [HN].

In the present work, in addition to $h$ on $\Theta$, we consider,
in particular, $h\big|_{\bar\Theta_{\rho}}$ and
$h\big|_{b\Theta_{\rho}}$, where
$$\eqalign{
&\bar\Theta_{\rho}=\Theta_{\rho}\cup b\Theta_{\rho},\cr
&\Theta_{\rho}=\{(k,l)\in\Theta:\ |Im\,k|=|Im\,l|<\rho\},\cr
&b\Theta_{\rho}=\{(k,l)\in\Theta:\ |Im\,k|=|Im\,l|=\rho\},\cr}\eqno(1.8)$$
where $\rho>0$. Note  that
$$dim\,\Theta=3d-4,\ \ dim\,b\Theta_{\rho}=3d-5.\eqno(1.9)$$
Using (1.4), (1.5) one can see that
$$h(k,l)\approx\hat v(p),\ \ p=k-l,\ \
(k,l)\in\Theta.\eqno(1.10)$$
in the Born approximation (that is in the linear approximation near zero
potential).
In addition, one can see also that
$$(k,l)\in\bar\Theta_{\rho}\Longrightarrow\ p=k-l\in \bar{\cal B}_{2\rho},
\eqno(1.11)$$
where
$$\eqalign{
&\bar{\cal B}_r={\cal B}_r\cup\pa {\cal B}_r,\cr
&{\cal B}_r=\{p\in\R^d:\ |p|<r\},\ \
\pa{\cal B}_r=\{p\in\R^d:\ |p|=r\},\ r>0.\cr}\eqno(1.12)$$

In the present work we consider, in particular, the following inverse
scattering problems for equation (1.1) under assumptions (1.2).

\vskip 2 mm
{\bf Problem 1.1.}

(a)\ Given $h$ on $\Theta$, find $v$ on $\R^d$;

(b)\ Given $h$ on $\bar\Theta_{\rho}$ for some (sufficiently great)
$\rho>0$, find $v$ on $\R^d$, at least, approximately;

(c)\ Given $h$ on $b\Theta_{\rho}$ for some (sufficiently great)
$\rho>0$, find $v$ on $\R^d$, at least, approximately.

Note that  Problems 1.1a, 1.1b make sense for any $d\ge 2$, whereas
Problem 1.1c is reasonable for $d\ge 3$ only:
$dim\,b\Theta_{\rho}< dim\,\R^d$ for $d=2$, see (1.9).

Note that: (1) any reconstruction method for Problem 1.1b with decaying
error as $\rho\to +\infty$ gives also a reconstruction method for
Problem 1.1a and (2) for $d\ge 3$, any reconstruction method for
Problem  1.1c gives also a reconstruction method for Problem 1.1b.

Note that in the Born approximation (1.10):
(a) Problem 1.1a is reduced to finding $v$ on $\R^d$ from
$\hat v$ on $\R^d$, (b)
Problem 1.1b is reduced to (approximate) finding $v$ on $\R^d$ from
$\hat v$ on ${\cal B}_{2\rho}$, (c)
Problem 1.1c for $d\ge 3$ is reduced to (approximate) finding $v$ on $\R^d$
from $\hat v$ on ${\cal B}_{2\rho}$,
where $\hat v$ is defined by (1.6).
Thus, in the Born approximation, Problem 1.1c for $d\ge 3$ (as well as
Problem 1.1b for $d\ge 2$) can be solved by the formula
$$\eqalign{
&v(x)=v_{appr}^{lin}(x,\rho)+v_{err}^{lin}(x,\rho),\cr
&v_{appr}^{lin}(x,\rho)=\int\limits_{{\cal B}_{2\rho}}e^{-ipx}\hat v(p)dp,\ \
v_{err}^{lin}(x,\rho)=\int\limits_{\R^d\b {\cal B}_{2\rho}}
e^{-ipx}\hat v(p)dp,\cr}\eqno(1.13)$$
where $x\in\R^d$. In addition, if, for example,
$$v\in W^{n,1}(\R^d)\ \ {\rm for\ some}\ \ n\in\N,\eqno(1.14)$$
and $\|v\|^{n,1}\le C$, where $W^{n,1}(\R^d)$ denotes the space of
n-times smooth functions on $\R^d$ in $L^1$-sense and $\|\cdot\|^{n,1}$
denotes some fixed standard norm in $W^{n,1}(\R^d)$, then
$$|\hat v(p)|\le c_1(n,d)C(1+|p|)^{-n},\ \ p\in\R^d,\eqno(1.15)$$
and, therefore, for $n>d$,
$$|v_{err}^{lin}(x,\rho)|\le c_2(n,d)C \rho^{-(n-d)},\ \ x\in\R^d,
\ \ \rho\ge 1,\eqno(1.16)$$
where $c_1(n,d)$, $c_2(n,d)$ are some fixed positive constants and
$v_{err}^{lin}(x,\rho)$ is the error term of (1.13).

Thus, in the Born approximation (1.10) (that is in the linear approximation
near zero potential) we have that:

(1) $h$ on $b\Theta_{\rho}$ for $d\ge 3$ (as well as $h$ on
$\bar\Theta_{\rho}$ for $d\ge 2$) stably determines $v_{appr}^{lin}(x,\rho)$ of
(1.13) and

(2) the error $v_{err}^{lin}(x,\rho)=v(x)-v_{appr}^{lin}(x,\rho)=
O(\rho^{-(n-d)})$ in the uniform norm as $\rho\to +\infty$ for n-times
smooth $v$ in the sense (1.14), where $n>d$. In particular,
$v_{err}^{lin}=O(\rho^{-\infty})$ in the uniform norm as $\rho\to +\infty$ for
 $v$  of the Schwartz class on $\R^d$.

The main results of the present work consist in global analogs for the
nonlinearized case for $d=3$ of the aforementioned Born-approximation
results for Problem 1.1c, see Theorem 2.1 and Corollary 2.1 of Section 2.
In particular, we give a stable approximate
solution of nonlinearized Problem 1.1c for $d=3$ and $v$ satisfying (1.14),
$n>d=3$, with the error term decaying as $O(\rho^{-(n-d)}\ln\rho)$ in the
uniform norm as $\rho\to +\infty$ (that is with almost the same decay
rate of the error for $\rho\to +\infty$ as in the linearized case near zero
potential, see (1.13), (1.16)).
The results of the present work were obtained in the framework of a
development of the $\bar\pa$-approach to inverse scattering at fixed energy
in dimension $d\ge 3$ of [BC], [HN], [No3], [No4],
with applications to the Gel'fand-Calderon inverse boundary value problem
via the reduction going back to [No1].
See Subsections 1.2, 1.3 and Sections 2, 3, 4, 5, 6 for details.

\vskip 2 mm
{\bf Remark 1.1.}
Note that if
$$\eqalign{
&v\in L^{\infty}(\R^d),\ \ ess\,\sup\limits_{x\in\R^d}
(1+|x|)^{d+\ep}|v(x)|\le C,\cr
&{\rm for\ some\ positive}\ \ \ep\ \ {\rm and}\ \ C,\cr}\eqno(1.17)$$
then (see [HN], [No1], [Na1], [No4]):
$$\eqalignno{
&\hat v(p)=\lim\limits_{\rho\to +\infty,\  k-l=p}h(k,l)\ \ {\rm for\ any}\ \
p\in\R^d,\ d\ge 3,&(1.18)\cr
&|\hat v(k-l)-h(k,l)|\le c_3(\ep,d)C^2\rho^{-1}\ \ {\rm as}\ \
\rho\to +\infty,&(1.19)\cr}$$
where $(k,l)\in b\Theta_{\rho}$, $c_3(\ep,d)$ is some positive constant.
Formulas (1.18), (1.19) show that the Born approximation (1.10) holds on
$b\Theta_{\rho}$ (and on $\Theta\b\Theta_{\rho}$) for any sufficiently
great $\rho$ (actually, for any sufficiently great $\rho$ in comparison with
$C$ of (1.17)). However, because of $O(\rho^{-1})$ in the right-hand side
of (1.19), formulas (1.18), (1.19) give no method to reconstruct $v$ on
$\R^d$ from $h$ on $b\Theta_{\rho}$ (or on $\bar\Theta_{\rho}$) with the
error
term decaying more rapidly than $O(\rho^{-1})$ in the uniform norm as
$\rho\to +\infty$ (even for $v$ of the Schwartz class on $\R^d$, $d\ge 3$).

\vskip 2 mm
{\bf Remark 1.2.}
On the other hand (in comparison with the result mentioned in Remark 1.1),
for sufficiently small potentials $v$, in [No4] we succeeded, in particular,
to give a stable method for solving Problem 1.1b for $d=3$ with the same
type rapid decay of the error term for $\rho\to +\infty$ as in formulas
(1.13), (1.16) for the linearized case near zero potential.
Moreover, in this result of [No4], $v$ is approximately reconstructed already
from non-overdetermined restriction $h\big|_{\Theta_{\rho}\cap\Gamma}$,
where $\Gamma\subset\Theta$, $dim\,\Gamma=d=3$. However, this result of
[No4] is local: the smallness of $v$ is used essentially.

\vskip 2 mm
{\bf Remark 1.3.}
We emphasize that before the present work no results were given, in general,
in the literature on solving Problems 1.1c and 1.1b for $d\ge 3$ with the
error term decaying more rapidly than $O(\rho^{-1})$ as $\rho\to +\infty$ even
 for $v$ of the Schwartz class on $\R^d$ (and even for the infinitely
smooth compactly supported case in the framework of sufficiently stable
rigorous algorithms). In addition, rapid decay of this error term is a
property of principal importance in the framework of applications of
methods for solving Problems 1.1 to the Gel'fand-Calderon inverse
boundary value problem (Problem 1.2) via the reduction of [No1], see the
next part of introduction.

Note that Problem 1.1a was considered for the first time in [BC] for $d=3$
from pure mathematical point of view without any physical applications.
No possibility to measure $h$ on $\Theta\b\{(0)\}$ directly in some
physical experiment is known at present
(here $\{0\}=\{(k,l)\in\Theta:\ |k|=|l|=0\}$).
However, as it was shown in [No1]
(see also [HN] (Note added in proof), [Na1], [No2]), Problems 1.1
naturally arise in the electrical impedance tomography and, more
generally, in connection with Problem 1.2 formulated in the next subsection.

\vskip 2 mm
{\it 1.2. The Gel'fand-Calderon problem.}
Consider the equation (1.1) in $D\subset\R^d$ only, where
$$\eqalign{
&D\ \ {\rm is\ an\ open\ bounded\ domain\ in}\ \ \R^d,\ \ d\ge 2,\cr
&{\rm with\ sufficiently\ regular\ boundary}\ \ \pa D,\cr
&v\ \ {\rm is\ a\ sufficiently\ regular\ function\ on}\ \ \bar D=D\cup\pa D.
\cr}\eqno(1.20)$$
For simplicity we assume also that
$$\eqalign{
&0\ \ {\rm is\ not\ a\ Dirichlet\ eigenvalue\ for} \cr
&{\rm the\ operator}\ \ -\Delta+v\ \ {\rm in}\ \ D.\cr}\eqno(1.21)$$
Consider the map $\Phi$ such that
$${\pa\psi\over \pa\nu}\big|_{\pa D}=\Phi\bigl(\psi\big|_{\pa D}\bigr)
\eqno(1.22)$$
for all sufficiently regular solutions $\psi$ of (1.1) in $\bar D$, where
$\nu$
is the outward normal to $\pa D$. The map $\Phi$ is called the
Dirichlet-to-Neumann map for equation (1.1) in $D$.
We consider the following  inverse
boundary value problem for equation (1.1) in $D$:

\vskip 2 mm
{\bf Problem 1.2.}
Given $\Phi$, find $v$.

This problem can be considered as the Gel'fand inverse boundary value problem
for the Schr\"odinger equation at zero energy (see [G], [No1]). This
problem can be also considered as a generalization of the Calderon problem
of the electrical impedance tomography
(see [C], [SU], [No1]).

One can see that
the Faddeev function $h$ of Problems 1.1 does not appear in Problem 1.2.
However, as it was shown in [No1] (see also [HN] (where this result of
[No1] was announced in Note added in proof),
[Na1], [No2]), if $h$ corresponds to equation (1.1) on $\R^d$, where
$v$ is the potential of Problem 1.2 on $D$ and $v\equiv 0$ on
$\R^d\b\bar D$,
then $h$ on $\Theta$ can be determined from the Dirichlet-to-Neumann map
$\Phi$ via the following formulas and equation:
$$\eqalignno{
&h(k,l)=(2\pi)^{-d}\int\limits_{\pa D}\int\limits_{\pa D}e^{-ilx}
(\Phi-\Phi_0)(x,y)\psi(y,k)dydx\ \ {\rm for}\ \ (k,l)\in\Theta,&(1.23)\cr
&\psi(x,k)=e^{ikx}+\int\limits_{\pa D}A(x,y,k)\psi(y,k)dy,\ \ x\in\pa D,
&(1.24)\cr
&A(x,y,k)=\int\limits_{\pa D}G(x-z,k)(\Phi-\Phi_0)(z,y)dz,\ \ x,y\in\pa D,
&(1.25)\cr
&G(x,k)=-(2\pi)^{-d}e^{ikx}\int\limits_{\R^d}{e^{i\xi x}d\xi\over
{\xi^2+2k\xi}},\ \ x\in\R^d,&(1.26)\cr}$$
where $k\in\C^d$, $k^2=0$ in (1.24)-(1.26), $\Phi_0$ denotes the
Dirichlet-to-Neumann map for equation (1.1) in $D$ with $v\equiv 0$, and
$(\Phi-\Phi_0)(x,y)$ is the Schwartz kernel of the integral operator
$\Phi-\Phi_0$. Note that (1.23), (1.25), (1.26) are explicit formulas,
whereas (1.24) is a linear integral equation (with parameter $k$) for
$\psi$ on $\pa D$. In addition, $G$ of (1.26) is the Faddeev's Green
function of [F1] for the Laplacian $\Delta$.
Formulas and equation (1.23)-(1.26) reduce Problem 1.2 to Problems 1.1.
In addition, from numerical point of view $h(k,l)$ for
$(k,l)\in\bar\Theta_{\rho}$ can be relatively easily
determined from $\Phi$
via (1.25), (1.24), (1.23) if $\rho$ is sufficiently small. However,
if $(k,l)\in\Theta\b\Theta_{\rho}$, where $\rho$ is sufficiently great,
then the  determination of $h(k,l)$ from $\Phi$ via (1.25), (1.24),
(1.23) is very unstable (especially on the step (1.24));
see, for example, [BRS], [No2], [No4]. This explains the principal importance
(mentioned in Remark 1.3) of the error term rapid decay as
$\rho\to +\infty$ in methods for solving Problems 1.1b and 1.1c.

\vskip 2 mm
{\it 1.3. Final remarks.}
In the present work we consider, mainly, Problems 1.1 and 1.2 for $d=3$.
The main results
of the present work are presented in Sections 2 and 6. Some of these results
were already mentioned above. Note that only restrictions in time prevent us
from generalizing all main results of the present work to the case
of the Schr\"odinger equation at arbitrary (not necessarily zero) fixed
energy $E$ for $d\ge 3$.

Note that results of the present work permit to complete (at least for
$d=3$) the proof of new stability estimates for Problem 1.2, $d\ge 3$,
announced as Theorem 2.2 of [NN]. We plan to return to this proof in a
separate article.

Our new global reconstruction for Problem 1.2 in 3D is summarized in
schemes (6.1), (6.2) of Section 6.
We expect that this reconstruction  can be implemented numerically in a 
similar way with
implementations developed in [ABR], where the parameter $\rho$ of (6.1),
(6.2) will play in some sense the role of the wave number $k_0$ of [ABR].

As regards results given in the literature on Problem 1.1, see [BC], [HN],
[GN], [Na1], [Na2], [No4] and references therein.

As regards results given in the literature on Problem 1.2
(in its Calderon or Gel'fand form), see [SU], [No1],
 [Al], [Na1], [Na2],  [Ma], [No2], [No4], [NN],  [HM], [Am]
and references therein.

\vskip 2 mm
{\bf 2.Main new results}

In the present work we consider, mainly, the three dimensional case $d=3$.
In addition, in the main considerations of the present work for $d=3$ our
basic assumption on $v$ consists in the following condition on its
Fourier transform:
$$\hat v\in L_{\mu}^{\infty}(\R^3)\cap {\cal C}(\R^3)\ \ {\rm for\ some\
real}\ \ \mu\ge 2,\eqno(2.1)$$
where $\hat v$ is defined by (1.6),
$$\eqalign{
&L_{\mu}^{\infty}(\R^d)=\{u\in L^{\infty}(\R^d):\ \|u\|_{\mu}< +\infty\},\cr
&\|u\|_{\mu}=ess \sup\limits_{p\in\R^d}(1+|p|)^{\mu}|u(p)|,\ \ \mu>0,\cr}
\eqno(2.2)$$
and ${\cal C}$ denotes the space of continuous functions. Actually, (2.1)
is a specification of (1.2).

Note that
$$\eqalign{
&v\in W^{n,1}(\R^d)\Longrightarrow\hat v\in L_{\mu}^{\infty}(\R^d)\cap
{\cal C}(\R^d),\cr
&\|\hat v\|_{\mu}\le c_4(n,d)\|v\|^{n,1}\ \ {\rm for}\ \ \mu=n,\cr}\eqno(2.3)
$$
where $W^{n,1}$, $L_{\mu}^{\infty}$ are the spaces of (1.14), (2.2).

Let
$$\eqalign{
&\Theta_{\rho,\tau}^{\infty}=\{(k,l)\in\Theta\b\bar\Theta_{\rho}:\
k-l\in {\cal B}_{2\rho\tau}\},\cr
&b\Theta_{\rho,\tau}=\{(k,l)\in b\Theta_{\rho}:\
k-l\in {\cal B}_{2\rho\tau}\},\cr}\eqno(2.4)$$
where $\Theta$, $\bar\Theta_{\rho}$,  $b\Theta_{\rho}$, ${\cal B}_r$ are
defined by (1.3), (1.8), (1.12), $\rho>0$, $0<\tau<1$. In (2.4) by
symbol $\infty$ we emphasize that  $\Theta_{\rho,\tau}^{\infty}$ is
unbounded with respect to $k,l$. One can see also that by definition
$$b\Theta_{\rho,\tau}\subset b\Theta_{\rho},\ \ \rho>0,\ 0<\tau<1.\eqno(2.5)$$

Our main new results on Problem 1.1 are summarized as Theorem 2.1 and
Corollary 2.1 below.

\vskip 2 mm
{\bf Theorem 2.1.}
{\it Let} $\hat v$ {\it satisfy} (2.1) {\it and} $\|\hat v\|_{\mu}\le C$.
{\it Let} $2\le\mu_0<\mu$, $0<\delta<1$,
$$0<\tau<\tau_1(\mu,\mu_0,C,\delta),\ \ \rho\ge\rho_1(\mu,\mu_0,C,\delta),
\eqno(2.6)$$
{\it where} $\tau_1$, $\rho_1$ {\it are the constants of} (4.36) {\it and, in
particular}, $0<\tau_1<1$. {\it Then} $h$ {\it on} $b\Theta_{\rho,\tau}$
{\it stably determines} $\hat v^{\pm}(\cdot,\tau,\rho)$ {\it via} (6.2)
({\it where the nonlinear integral equation} (4.31) {\it is solved by
successive approximations) and}
$$|\hat v(p)-\hat v^{\pm}(p,\tau,\rho)|\le c_5(\mu,\mu_0,\tau,\delta)C^2
(1+|p|)^{-\mu_0}\rho^{-(\mu-\mu_0)}\ \ {\it for}\ \ p\in {\cal B}_{2\tau\rho},
\eqno(2.7)$$
{\it where} $c_5$ {\it is the constant of} (6.4), see Section 6.

Construction (6.2) is actually a definition of
$\hat v^{\pm}(\cdot,\tau,\rho)$. We consider $\hat v^{\pm}(\cdot,\tau,\rho)$
as an approximation to $\hat v$ on  ${\cal B}_{2\tau\rho}$. The error
between  $\hat v$ and $\hat v^{\pm}(\cdot,\tau,\rho)$ on
${\cal B}_{2\tau\rho}$ is estimated in (2.7). Estimate (2.7) is especially
interesting for $\rho\to +\infty$ at fixed $\tau$.

Theorem 2.1 follows from Proposition 4.2, Corollary 4.1 and formulas (5.2),
(5.3), (5.6) (and related results of Sections 4 and 5).
A more detailed   version of (2.7) is given by (6.3). The stability
mentioned in Theorem 2.1 follows from estimate (6.10).
See Sections 4, 5, 6 for additional details.

An  outline of the reconstruction
formalized in Theorem 2.1 is given also at the end of the present section.

Let
$$v^{\pm}(x,\tau,\rho)=\int\limits_{{\cal B}_{2\tau\rho}}
e^{-ipx}\hat v^{\pm}(p,\tau,\rho)dp,\ x\in\R^3,\eqno(2.8)$$
where $\hat v^{\pm}(p,\tau,\rho)$ is the approximation of Theorem 2.1.

Formula (2.3) and Theorem 2.1 imply the following

\vskip 2 mm
{\bf Corollary 2.1.}
{\it Let} $v$ {\it satisfy} (1.14), $n>d=3$, {\it and} $\|v\|^{n,1}\le D$.
{\it Let} $\tau$ {\it and} $\rho$ {\it satisfy} (2.6) {\it for}
$\mu=n$, $\mu_0=3$, $C=c_4(n,3)D$ {\it and fixed} $\delta\in ]0,1[$.
{\it Then} $h$ {\it on} $b\Theta_{\rho,\tau}$ {\it stably determines}
$v^{\pm}(\cdot,\tau,\rho)$ {\it via} (6.2), (2.8) {\it and}
$$|v(x)-v^{\pm}(x,\tau,\rho)|\le (c_6(n,\tau)D+c_7(n,\tau,\delta)D^2
\ln\,(1+2\tau\rho)) \rho^{-(n-3)}\ \ {\it for}\ \ x\in\R^3,\eqno(2.8)$$
{\it where the constants} $c_6$, $c_7$ {\it are simply related with} $c_2$
{\it and} $c_5$. {\it In addition, in particular},
$$\eqalign{
&\|v-v^{\pm}(\cdot,\tau,\rho)\|_{L^{\infty}(\R^3)}=O(\rho^{-(n-3)}\ln\,\rho)\
 \ {\it for}\ \ \rho\to +\infty\ \ {\it for\ fixed}\ \ \tau,\cr
&0<\tau<\tau_1(n,3,c_4(n,3)D,\delta).\cr}\eqno(2.9)$$

We consider Theorem 2.1 as a global nonlinear analog of the result that in
the Born approximation $h\big|_{b\Theta_{\rho}}$ is reduced to a $\hat v$ on
${\cal B}_{2\rho}$ for $d\ge 3$.  One can see, in particular, that in
Theorem 2.1 we even do not try to reconstruct $\hat v$ on
$\R^3\b {\cal B}_{2\rho}$ from  $h$ on $b\Theta_{\rho}$.

We consider Corollary 2.1 as a global nonlinear analog of the Born
approximation result (for Problem 1.1c) consisting in formulas (1.13), (1.16).

In the derivations of the present work we rewrite $h$ on $\Theta$,
$\bar\Theta_{\rho}$, $b\Theta_{\rho}$,
$\Theta_{\rho,\tau}^{\infty}$ and $b\Theta_{\rho,\tau}$
as $H$ on $\Omega$, $\bar\Omega_{\rho}$, $\Omega_{\rho,\tau}^{\infty}$
and
$b\Omega_{\rho,\tau}$ (respectively), where $h$ is related with $H$ by (1.4),
$$\eqalign{
&\Omega=\{k\in\C^d,\ p\in\R^d:\ k^2=0,\ p^2=2kp\},\cr
&\bar\Omega_{\rho}=\Omega_{\rho}\cup b\Omega_{\rho},\cr
&\Omega_{\rho}=\{(k,p)\in\Omega:\ |Im\,k|<\rho\},\cr
&b\Omega_{\rho}=\{(k,p)\in\Omega:\ |Im\,k|=\rho\},\cr
&\Omega_{\rho,\tau}^{\infty}=\{(k,p)\in\Omega\b\bar\Omega_{\rho}:\
p\in {\cal B}_{2\rho\tau}\},\cr
&b\Omega_{\rho,\tau}=\{(k,p)\in b\Omega_{\rho}:\
p\in {\cal B}_{2\rho\tau}\},\cr}\eqno(2.10)$$
where $\rho>0$, $0<\tau<1$.

Note that
$$\eqalign{
&\Omega\approx\Theta,\ \bar\Omega_{\rho}\approx\bar\Theta_{\rho},\
b\Omega_{\rho}\approx b\Theta_{\rho},\cr
&\Omega_{\rho,\tau}^{\infty}\approx\Theta_{\rho,\tau}^{\infty},\ \
b\Omega_{\rho,\tau}\approx b\Theta_{\rho,\tau}.\cr}\eqno(2.11)$$
or more precisely
$$\eqalign{
&(k,p)\in\Omega\Longrightarrow (k,k-p)\in\Theta,\ \
(k,l)\in\Theta\Longrightarrow (k,k-l)\in\Omega\cr
&{\rm and\ the\ same\ for}\ \
\bar\Omega_{\rho},\ b\Omega_{\rho},\ \Omega_{\rho,\tau}^{\infty},\
b\Omega_{\rho,\tau}\ \
{\rm and}\ \
\bar\Theta_{\rho},\ b\Theta_{\rho},\ \Theta_{\rho,\tau}^{\infty},\
b\Theta_{\rho,\tau},\cr
&{\rm respectively},\ \
{\rm in\ place\ of}\ \Omega\ \ {\rm and}\ \ \Theta.\cr}
\eqno(2.12)$$

An outline of the reconstruction formalized in Theorem 2.1 consists in the
following:
\item{1.} We rewrite $h$ on $b\Theta_{\rho,\tau}\subset b\Theta_{\rho}$ as
$H$ on $b\Omega_{\rho,\tau}\subset b\Omega_{\rho}$
as mentioned above.
\item{2.} We consider $H$  on $b\Omega_{\rho}$ as boundary data for $H$ on
$\Omega\b\Omega_{\rho}$, which solves the non-linear $\bar\pa$-equation (3.5)
with estimates (3.2), (3.3).
\item{3.} Using this $\bar\pa$-equation and these estimates we obtain the
non-linear integral equation (4.31) for finding $\tilde H_{\rho,\tau}$ on
$\Omega_{\rho,\tau}^{\infty}$ from $H$ on $b\Omega_{\rho,\tau}$, where
$\tilde H_{\rho,\tau}$  approximates $H$ on $\Omega_{\rho,\tau}$
(with estimate (4.31)).
\item{4.} The function $\tilde H_{\rho,\tau}$ determines
$\hat v^{\pm}(\cdot,\tau,\rho)$ by formulas (4.32), where
$\hat v^{\pm}(\cdot,\tau,\rho)$   approximates $\hat v$ on
${\cal B}_{2\tau\rho}$.

As it was already mentioned in introduction, the results of the present work
were obtained in the framework of a development of the $\bar\pa$-approach
to inverse scattering at fixed energy in dimension $d\ge 3$ of [BC], [HN],
[No3], [No4]. In particular, there is a considerable similarity between the
reconstruction scheme of the present work for the case of Problem 1.1c
 and the reconstruction scheme of [No3] for the case of inverse scattering
at fixed positive energy $E$. Actually, in the present work the parameter
$\rho$ of Problem 1.1 plays the role of $E^{1/2}$ of [No3]. Some results
of [BC], [HN], [No4] we use in the present work are recalled in the next
section.

\vskip 2 mm
{\bf 3. Background results}

For simplicity always in this section we assume that $d=3$.

\vskip 2 mm
{\it 3.1. Estimates for $H$ on $\Omega\b\Omega_{\rho}$}.
Let $v$ satisfy (2.1) and $\|\hat v\|_{\mu}\le C$. Let
$$\eta(C,\rho,\mu)\buildrel \rm def \over = a(\mu)C(\ln\,\rho)^2\rho^{-1}<1,\
\ \ln\,\rho\ge 2,\eqno(3.1)$$
where $a(\mu)$ is the constant $c_2(\mu)$ of [No4]. Then (according to
[No4]):
$$\eqalignno{
&H\in {\cal C}(\Omega\b\Omega_{\rho}),&(3.2)\cr
&|H(k,p)|\le {C\over (1-\eta(C,\rho,\mu))(1+|p|)^{\mu}},\ \
(k,p)\in\Omega\b\Omega_{\rho},&(3.3)\cr
&\hat v(p)=\lim\limits_{|k|\to\infty,\ (k,p)\in\Omega}H(k,p),\ \ p\in\R^3,
&(3.4)\cr}$$
where $|k|=((Re\,k)^2+(Im\,k)^2)^{1/2}$. (These and some additional
estimates on $H$ are given in Proposition 3.2 of [No4].) Note that, for
sufficiently regular $v$ on $\R^3$ with sufficient decay at infinity,
formula (3.4) was obtained for the first time in [HN].

\vskip 2 mm
{\it 3.2. The $\bar\pa$-equation for $H$ on $\Omega\b\Omega_{\rho}$.}
Let $v$ satisfy (2.1), $\|\hat v_{\mu}\|\le C$, and (3.1) hold. Then
(see [No4]):
$$\eqalign{
&\bar\pa_kH(k,p)\big|_{\Omega\b\Omega_{\rho}}=\cr
&\sum_{j=1}^3\biggl(-2\pi\int\limits_{\xi\in S_k}\xi_jH(k,-\xi)H(k+\xi,p+\xi)
{ds\over |Im\,k|^2}\biggr)d\bar k_j\big|_{\Omega\b\Omega_{\rho}},\cr}
\eqno(3.5)$$
where
$$S_k=\{\xi\in\R^3:\ \xi^2+2k\xi=0\},\eqno(3.6)$$
$ds$ is arc-length measure on the circle $S_k$ in $\R^3$.
Actually, under some stronger assumptions on $v$ than in the present
subsection, the $\bar\pa$- equation (3.5) was obtained for the first time
in [BC].

An important property of the $\bar\pa$-equation (3.5) is that (3.5) can be
considered for $H$ on $\Omega\b\Omega_{\rho}$ only, see, in particular,
formulas (3.24) of Subsection 3.4.

\vskip 2 mm
{\it 3.3. Coordinates on $\Omega\b\Omega_{\rho}$.}
Let
$$\Omega_{\nu}=\{(k,p)\in\Omega:\ p\not\in {\cal L}_{\nu}\},\eqno(3.7)$$
where
$${\cal L}_{\nu}=\{p\in\R^3:\ p=t\nu,\ t\in\R\},\ \ \nu\in\S^2.\eqno(3.8)$$
For $p\in\R^3\b {\cal L}_{\nu}$ we consider $\theta(p)$ and $\omega(p)$ such
that
$$\eqalign{
&\theta(p),\omega(p)\ \ {\rm smoothly\ depend\ on}\ \ p\in\R^3\b
{\cal L}_{\nu},\cr
&{\rm take\ values\ in}\ \ \S^2,\ \ {\rm and}\cr
&\theta(p)p=0,\ \omega(p)p=0,\ \theta(p)\omega(p)=0.\cr}\eqno(3.9)$$
Assumptions (3.9) imply that
$$\omega(p)={p\times\theta(p)\over |p|}\ \ {\rm for}\ \
p\in\R^3\b {\cal L}_{\nu} \eqno(3.10a)$$
or
$$\omega(p)=-{p\times\theta(p)\over |p|}\ \ {\rm for}\ \
p\in\R^3\b {\cal L}_{\nu}, \eqno(3.10b)$$
where $\times$ denotes vector product.

To satisfy (3.9), (3.10a) we can take
$$\theta(p)={\nu\times p\over |\nu\times p|},\ \omega(p)=
{p\times\theta(p)\over |p|},\ p\in\R^3\b {\cal L}_{\nu}.\eqno(3.11)$$

Let $\theta,\omega$ satisfy (3.9). Then (according to [No4]) the following
formulas give a diffeomorphism between $\Omega_{\nu}$ and
$(\C\b 0)\times (\R^3\b {\cal L}_{\nu})$:
$$(k,p)\to (\lambda,p),\ \ {\rm where}\ \
\lambda=\lambda(k,p)={2k(\theta(p)+i\omega(p))\over i|p|},\eqno(3.12a)$$
$$\eqalign{
&(\lambda,p)\to (k,p),\ \ {\rm where}\ \ k=k(\lambda,p)=\kappa_1(\lambda,p)
\theta(p)+\kappa_2(\lambda,p)\omega(p)+{p\over 2},\cr
&\kappa_1(\lambda,p)={i|p|\over 4}(\lambda+{1\over \lambda}),
\ \  \kappa_2(\lambda,p)={|p|\over 4}(\lambda-{1\over \lambda}),\cr}
\eqno(3.12b)$$
where $(k,p)\in\Omega_{\nu}$,
$(\lambda,p)\in (\C\b 0)\times (\R^3\b {\cal L}_{\nu})$.
In addition, formulas (3.12a), (3.12b) for $\lambda(k)$ and $k(\lambda)$
at fixed $p\in\R^3\b {\cal L}_{\nu}$
give a diffeomorphism between $Z_p=\{k\in\C^3:\ (k,p)\in\Omega\}$ for fixed
$p$ and $\C\b 0$.

In addition, for $k$ and $\lambda$ of (3.12) we have that
$$|Im\,k|={|p|\over 4}\bigl(|\lambda|+{1\over |\lambda|}\bigr),\ \
|Re\,k|={|p|\over 4}\bigl(|\lambda|+{1\over |\lambda|}\bigr),\eqno(3.13)$$
where $(k,p)\in\Omega_{\nu}$,
$(\lambda,p)\in (\C\b 0)\times (\R^3\b {\cal L}_{\nu})$.

Let
$$\Omega_{\rho,\tau,\nu}^{\infty}=\Omega_{\rho,\tau}^{\infty}\cap\Omega_{\nu},
\eqno(3.14)$$
where $\Omega_{\rho,\tau}^{\infty}$, $\Omega_{\nu}$ are defined as in (2.10),
(3.7). Let
$$\Lambda_{\rho,\nu}=\{(\lambda,p):\ \lambda\in {\cal D}_{\rho/|p|},\ \
p\in\R^3\b {\cal L}_{\nu}\},\eqno(3.15)$$
$$\eqalign{
&\Lambda_{\rho,\tau,\nu}=\{(\lambda,p):\ \lambda\in {\cal D}_{\rho/|p|},\ \
p\in\R^3\b {\cal L}_{\nu},\ \ |p|<2\tau\rho\},\cr
&b\Lambda_{\rho,\tau,\nu}=\{(\lambda,p):\ \lambda\in {\cal T}_{\rho/|p|},\ \
p\in\R^3\b {\cal L}_{\nu},\ \ |p|<2\tau\rho\},\cr}$$
where $\rho>0$, $0<\tau<1$, $\nu\in\S^2$,
$${\cal D}_r=\{\lambda\in\C\b 0:\ {1\over 4}(|\lambda|+|\lambda|^{-1})>r\},\
r>0.\eqno(3.16)$$
$${\cal T}_r=\{\lambda\in\C:\ {1\over 4}(|\lambda|+|\lambda|^{-1})=r\},\
r\ge 1/2.\eqno(3.17)$$
Using (3.13) one can see that  formulas (3.12) give also the following
diffeomorphisms
$$\eqalign{
&\Omega_{\nu}\b\bar\Omega_{\rho}\approx\Lambda_{\rho,\nu},\ \
\Omega_{\rho,\tau,\nu}^{\infty}\approx\Lambda_{\rho,\tau,\nu},\cr
&b\Omega_{\rho,\tau}\cap\Omega_{\nu}\approx b\Lambda_{\rho,\tau,\nu},\cr
&Z_{p,\rho}^{\infty}=\{k\in\C^3:\ (k,p)\in\Omega_{\nu}\b\bar\Omega_{\rho}\}
\approx {\cal D}_{\rho/|p|}\ \ {\rm for\ fixed}\ \ p,\cr}\eqno(3.18)$$
where $\rho>0$, $0<\tau<1$, $\nu\in\S^2$.

In [No4] $\lambda,p$ of (3.12) were used as coordinates on
$\Omega$. In the present work we use them also as coordinates on
$\Omega\b\Omega_{\rho}$ (or more precisely on $\Omega_{\nu}\b\Omega_{\rho}$).

\vskip 2 mm
{\it 3.4. $\bar\pa$-equation for $H$ in the $\lambda,p$ coordinates and
some related estimate.}
Let $\lambda,p$ be the coordinates of Subsection 3.3, where $\theta$,
$\omega$ satisfy (3.9), (3.10). Then (see Lemma 5.1 of [No4]) in these
coordinates the $\bar\pa$-equation (3.5) for $p\ne 0$ takes the form:
$$\eqalign{
&{\pa\over \pa\bar\lambda}H(k(\lambda,p),p)=-{\pi\over 4}
\int_{-\pi}^{\pi}\biggl({|p|\over 2} {(|\lambda|^2-1)\over \bar\lambda
|\lambda|}(\cos\v-1)-|p|{1\over \bar\lambda}\sin\v\biggr)\times\cr
&H(k(\lambda,p),-\xi(\lambda,p,\v))H(k(\lambda,p,\v),p+\xi(\lambda,p,\v))
d\v\cr}\eqno(3.19)$$
for $(\lambda,p)\in\Lambda_{\rho,\nu}$, where
$k(\lambda,p)$ is defined in (3.12b) (and also depends on $\nu$, $\theta$,
$\omega$),  $\Lambda_{\rho,\nu}$ is defined in (3.15),
$$\eqalignno{
&\xi(\lambda,p,\v)=Re\,k(\lambda,p)(\cos\v-1)+k^{\perp}(\lambda,p)\sin\v,
&(3.20)\cr
&k^{\perp}(\lambda,p)={Im\,k(\lambda,p)\times Re\,k(\lambda,p)\over
|Im\,k(\lambda,p)|},&(3.21)\cr}$$
where $\times$ in (3.21) denotes vector product.

Note that (3.19) can be written as
$${\pa\over \pa\bar\lambda}H(k(\lambda,p),p)=\{H,H\}(\lambda,p),\ \
(\lambda,p)\in\Lambda_{\rho,\nu},\eqno(3.22)$$
where
$$\eqalign{
&\{U_1,U_2\}(\lambda,p)=-{\pi\over 4}
\int_{-\pi}^{\pi}\biggl({|p|\over 2} {{|\lambda|^2-1}\over \bar\lambda
|\lambda|}(\cos\v-1)-{|p|\over \bar\lambda}\sin\v\biggr)\times\cr
&U_1(k(\lambda,p),-\xi(\lambda,p,\v))U_2(k(\lambda,p)+\xi(\lambda,p,\v),
p+\xi(\lambda,p,\v))d\v,\cr}\eqno(3.23)$$
where $U_1$, $U_2$ are test functions on $\Omega\b\bar\Omega_{\rho}$,
$k(\lambda,p)$, $\xi(\lambda,p,\v)$ are defined by (3.12b), (3.20),
$(\lambda,p)\in\Lambda_{\rho,\nu}$.
Note that in the left-hand side of (3.19), (3.22)
$$(k(\lambda,p),p)\in\Omega_{\nu}\b\bar\Omega_{\rho} \eqno(3.24a)$$
and in the right-hand side of (3.19), (3.23)
$$\eqalign{
&(k(\lambda,p),-\xi(\lambda,p,\v))\in\Omega\b\bar\Omega_{\rho},\cr
&(k(\lambda,p)+\xi(\lambda,p,\v),p+\xi(\lambda,p,\v))
\in\Omega\b\bar\Omega_{\rho},\cr}\eqno(3.24b)$$
where $(\lambda,p)\in\Lambda_{\rho,\nu}$, $\v\in [-\pi,\pi]$.

Let $U_1,U_2\in L_{\mu}^{\infty}(\Omega\b\bar\Omega_{\rho})$,
$\mu\ge 2$, where
$$\eqalign{
&L^{\infty}_{\mu}(\Omega\b\bar\Omega_{\rho})=
\{U\in L^{\infty}(\Omega\b\bar\Omega_{\rho}):\ |||U|||_{\rho,\mu}<+\infty
\},\cr
&|||U|||_{\rho,\mu}=ess\,\sup_{(k,p)\in\Omega\b\bar\Omega_{\rho}}(1+|p|)^{\mu}
|U(k,p)|,\ \mu>0.\cr}\eqno(3.25)$$

Let $\{U_1,U_2\}$ be defined by (3.23). Then (as a corollary of
Lemma 5.2 of [No4]):
$$\{U_1,U_2\}\in L^{\infty}_{local}(\Lambda_{\rho,\nu}) \eqno(3.26)$$
and
$$\eqalign{
&|\{U_1,U_2\}(\lambda,p)|\le {|||U_1|||_{\rho,\mu}|||U_2|||_{\rho,\mu}\over
(1+|p|)^{\mu}}b(\mu,|\lambda|,|p|)\cr
&{\rm for\ almost\ all}\ \ (\lambda,p)\in\Lambda_{\rho,\nu},\ \
b(\mu,|\lambda|,|p|)=\cr
&\biggl({b_1(\mu)|\lambda|\over (|\lambda|^2+1)^2}+
{b_2(\mu)|p|||\lambda|^2-1|\over |\lambda|^2(1+|p|(|\lambda|+|\lambda|^{-1}))
^2}+
{b_3(\mu)|p|\over |\lambda|(1+|p|(|\lambda|+|\lambda|^{-1}))}\biggr),\cr}
\eqno(3.27)$$
where $\Lambda_{\rho,\nu}$ is defined in (3.15), $b_1(\mu)$, $b_2(\mu)$,
$b_3(\mu)$ are the constants  $c_3(\mu)$, $c_4(\mu)$, $c_5(\mu)$ of [No4].

\vskip 2 mm
{\bf 4. Approximate finding $H$ on $\Omega_{\rho,\tau}^{\infty}$ from $H$ on
$b\Omega_{\rho,\tau}$}

We recall that $\Omega_{\rho,\tau}^{\infty}$ and  $b\Omega_{\rho,\tau}$ were
defined in Section 2, see formulas (2.10).  We assume that $d=3$.

Consider $\chi_rH$, where $\chi_r$ denotes the multiplication operator by
the function
$$\chi_r(p)=1\ \ {\rm for}\ \ |p|<r,\ \chi_r(p)=0\ \ {\rm for}\ \
|p|\ge r,\ \ {\rm where}\ \ p\in\R^3,\ r>0.\eqno(4.1)$$
Note that
$$\eqalign{
&\chi_{2\tau\rho}H(k,p)=H(k,p)\ \ {\rm for}\ \ (k,p)\in
\Omega_{\rho,\tau}^{\infty}\cr
&\chi_{2\tau\rho}H(k,p)=0\ \ {\rm for}\ \ (k,p)\in
(\Omega\b\bar\Omega_{\rho})\b\Omega_{\rho,\tau}^{\infty},\cr}\eqno(4.2)$$
where $\rho>0$, $\tau\in ]0,1[$.

As a corollary of (3.5), (3.22), (3.23) we have that
$${\pa\over \pa\bar\lambda}\chi_{2\tau\rho}H(k(\lambda,p),p)=
\{\chi_{2\tau\rho}H,\chi_{2\tau\rho}H\}(\lambda,p)+R_{\rho,\tau}(\lambda,p),
\eqno(4.3)$$
$$\eqalign{
&R_{\rho,\tau}(\lambda,p)=\cr
&\{(1-\chi_{2\tau\rho})H,\chi_{2\tau\rho}H\}(\lambda,p)+
\{\chi_{2\tau\rho}H,(1-\chi_{2\tau\rho})H\}(\lambda,p)+\cr
&\{(1-\chi_{2\tau\rho})H,(1-\chi_{2\tau\rho})H\}(\lambda,p)\cr}\eqno(4.4)$$
for $(\lambda,p)\in\Lambda_{\rho,\tau,\nu}$ of (3.15).

Because of the remainder $R_{\rho,\tau}$ of   (4.3), (4.4), the
$\bar\pa$-equation (3.5), (4.3) is only an approximate $\bar\pa$- equation
for $\chi_{2\tau\rho}H=H$ on $\Omega_{\rho,\tau}^{\infty}$ or on
$\Lambda_{\rho,\tau,\nu}$ in the coordinates $\lambda$, $p$. However,
$R_{\rho,\tau}$ rapidly vanishes when $\rho$ increases for fixed
$\tau\in ]0,1[$; see Lemma 4.1.

For approximate finding $H$ on $\Omega_{\rho,\tau}^{\infty}$ from $H$ on
$\Omega_{\rho,\tau}^{\infty}$ we proceed from (3.2), (3.3), (4.3), (4.4),
(3.26), (3.27) and the following formulas
$$\eqalignno{
&u_+(\lambda)={1\over 2\pi i}\int\limits_{{\cal T}_r^+}u_+(\zeta)
{d\zeta\over {\zeta-\lambda}}-{1\over \pi}
\int\!\!\!\int\limits_{{\cal D}_r^+}
{\pa u_+(\zeta)\over \pa\bar\zeta}
{d\,Re\,\zeta\ d\,Im\,\zeta\over {\zeta-\lambda}},\ \
\lambda\in {\cal D}_r^+,&(4.5a)\cr
&u_-(\lambda)=-{1\over 2\pi i}\int\limits_{{\cal T}_r^-}u_-(\zeta)
{\lambda d\zeta\over \zeta(\zeta-\lambda)}-
{1\over \pi}\int\!\!\!\int\limits_{{\cal D}_r^-}
{\pa u_-(\zeta)\over \pa\bar\zeta}{\lambda d\,Re\,\zeta\ d\,Im\,\zeta
\over \zeta(\zeta-\lambda)},\ \ \lambda\in {\cal D}_r^-,&(4.5b)\cr}$$
where
$$\eqalign{
&{\cal D}_r^{\pm}=\{\lambda\in\C\b 0:\ {1\over 4}(|\lambda|+|\lambda|^{-1})>r,
\ \ |\lambda|^{\pm 1}<1\},\cr
&{\cal T}_r^{\pm}=\{\lambda\in\C:\ {1\over 4}(|\lambda|+|\lambda|^{-1})=r,
\ \ |\lambda|^{\pm 1}\le 1\},\ \ r>1/2,\cr}\eqno(4.6)$$
$u_+(\lambda)$ is continuous and bounded on ${\cal D}_r^+\cup {\cal T}_r^+$,
$\pa u_+(\lambda)/\pa\bar\lambda$ is bounded on ${\cal D}_r^+$,
$u_-(\lambda)$ is continuous and bounded on
${\cal D}_r^-\cup {\cal T}_r^-$,
$\pa u_-(\lambda)/\pa\bar\lambda$ is bounded on ${\cal D}_r^-$, and
$\pa u_-(\lambda)/\pa\bar\lambda=O(|\lambda|^{-2})$ as $|\lambda|\to\infty$
(and where the integrals along ${\cal T}_r^{\pm}$ are taken in the
counter-clockwise direction). The aforementioned assumptions on
$u_{\pm}$ in (4.5) can be somewhat weakened. Formulas (4.5) follow from the
well-known Cauchy-Green formula
$$u(\lambda)={1\over 2\pi i}\int\limits_{\pa {\cal D}}u(\zeta)
{d\zeta\over {\zeta-\lambda}}
-{1\over \pi}\int\!\!\!\int\limits_{\cal D}
{\pa u(\zeta)\over \pa\bar\zeta}
{d\,Re\,\zeta\ d\,Im\,\zeta\over {\zeta-\lambda}},\ \
\lambda\in {\cal D},\eqno(4.7)$$
where $\cal D$ is a bounded open domain in $\C$ with sufficiently regular
boundary and $u$ is a sufficiently regular function on
$\bar {\cal D}={\cal D}\cup\pa {\cal D}$.

Let
$$H(\lambda,p)=H(k(\lambda,p),p),\ \ (\lambda,p)\in (\C\b 0)\times
(\R^3\b {\cal L}_{\nu}),\eqno(4.8)$$
where $\lambda,p$ are the coordinates of Subsection 3.3 under assumption
(3.10a).

Let
$$\eqalign{
&\Lambda^{\pm}_{\rho,\tau,\nu}=\{(\lambda,p):\
\lambda\in {\cal D}^{\pm}_{\rho/|p|},\
p\in {\cal B}_{2\tau\rho}\b {\cal L}_{\nu}\},\cr
&b\Lambda^{\pm}_{\rho,\tau,\nu}=\{(\lambda,p):\
\lambda\in {\cal T}^{\pm}_{\rho/|p|},\
p\in {\cal B}_{2\tau\rho}\b {\cal L}_{\nu}\},\cr}\eqno(4.9)$$
where  ${\cal B}_r$,  ${\cal L}_{\nu}$, ${\cal D}_r^{\pm}$, ${\cal T}_r^{\pm}$
are defined by
(1.12) for $d=3$, (3.8), (4.6), $\rho>0$, $\tau\in ]0,1[$, $\nu\in\S^2$.

Note that
$$\Lambda_{\rho,\tau,\nu}=
\Lambda_{\rho,\tau,\nu}^+\cup\Lambda_{\rho,\tau,\nu}^-,\ \
\Lambda_{\rho,\tau,\nu}^+\cap\Lambda_{\rho,\tau,\nu}^-=\emptyset,\ \
b\Lambda_{\rho,\tau,\nu}=
b\Lambda_{\rho,\tau,\nu}^+\cup\Lambda_{\rho,\tau,\nu}^-,\eqno(4.10)$$
where $\Lambda_{\rho,\tau,\nu}$, $b\Lambda_{\rho,\tau,\nu}$
were defined in (3.15), $\rho>0$, $\tau\in ]0,1[$, $\nu\in\S^2$.

As a corollary of (3.2), (3.3), (4.3), (4.4), (3.26), (3.27), (4.5), we
obtain the following

\vskip 2 mm
{\bf Proposition 4.1.}
{\it Let} $v$ {\it and} $\rho$ {\it satisfy the same assumptions that in
Subsection  3.1. Let}  $H(\lambda,p)$ {\it be defined by} (4.8). {\it Then}
$H=H(\lambda,p)$ {\it as a function of} $(\lambda,p)\in
\Lambda_{\rho,\tau,\nu}$ {\it of} (4.10), {\it where} $\tau\in ]0,1[$,
{\it satisfies the following nonlinear integral equation}
$$H=H^0+M_{\rho,\tau}(H)+Q_{\rho,\tau},\ \ \tau\in ]0,1[, \eqno(4.11)$$
{\it where}
$$\eqalignno{
&H^0(\lambda,p)={1\over 2\pi i}\int\limits_{{\cal T}^+_{\rho/|p|}}
H(\zeta,p){d\zeta\over {\zeta-\lambda}},\
(\lambda,p)\in\Lambda^+_{\rho,\tau,\nu},&(4.12a)\cr
&H^0(\lambda,p)=-{1\over 2\pi i}\int\limits_{{\cal T}^-_{\rho/|p|}}
H(\zeta,p){\lambda d\zeta\over \zeta(\zeta-\lambda)},\
(\lambda,p)\in\Lambda^-_{\rho,\tau,\nu},
&(4.12b)\cr}$$
{\it where} ${\cal T}^{\pm}_r$ {\it are defined by} (4.6);
$$\eqalign{
&M_{\rho,\tau}(U)(\lambda,p)=M_{\rho,\tau}^+(U)(\lambda,p)=\cr
&-{1\over \pi}\int\!\!\!\int\limits_{{\cal D}_{\rho/|p|}^+}
(U,U)_{\rho,\tau}(\zeta,p)
{d\,Re\,\zeta\ d\,Im\,\zeta\over {\zeta-\lambda}},\ \
(\lambda,p)\in\Lambda^+_{\rho,\tau,\nu},\cr}\eqno(4.13a)$$
$$\eqalign{
&M_{\rho,\tau}(U)(\lambda,p)=M_{\rho,\tau}^-(U)(\lambda,p)=\cr
&-{1\over \pi}\int\!\!\!\int\limits_{{\cal D}_{\rho/|p|}^-}
(U,U)_{\rho,\tau}(\zeta,p)
{\lambda d\,Re\,\zeta\ d\,Im\,\zeta\over \zeta(\zeta-\lambda)},\ \
(\lambda,p)\in\Lambda^-_{\rho,\tau,\nu},\cr}\eqno(4.13b)$$
$$\eqalign{
&(U_1,U_2)_{\rho,\tau}(\zeta,p)=
\{\chi_{2\tau\rho}U_1^{\prime},\chi_{2\tau\rho}U_2^{\prime}\}(\zeta,p),\
(\zeta,p)\in\Lambda_{\rho,\tau,\nu},\cr
&\chi_{2\tau\rho}U_j^{\prime}(k,p)=U_j(\lambda(k,p),p),\
(k,p)\in\Omega_{\rho,\tau,\nu}^{\infty},\cr
&\chi_{2\tau\rho}U_j^{\prime}(k,p)=0,\ |p|\ge 2\tau\rho,\ j=1,2,\cr}
\eqno(4.14)$$
{\it where} $U,U_1,U_2$ {\it are test functions on} $\Lambda_{\rho,\tau,\nu}$,
 $\{\cdot,\cdot\}$ {\it is defined by} (3.23), $\lambda(k,p)$ {\it is
defined in} (3.12a);
$$\eqalignno{
&Q_{\rho,\tau}(\lambda,p)=
-{1\over \pi}\int\!\!\!\int\limits_{{\cal D}_{\rho/|p|}^+}
R_{\rho,\tau}(\zeta,p)
{d\,Re\,\zeta\ d\,Im\,\zeta\over {\zeta-\lambda}},\ \
(\lambda,p)\in\Lambda^+_{\rho,\tau,\nu},&(4.15a)\cr
&Q_{\rho,\tau}(\lambda,p)=
-{1\over \pi}\int\!\!\!\int\limits_{{\cal D}_{\rho/|p|}^-}
R_{\rho,\tau}(\zeta,p)
{\lambda d\,Re\,\zeta\ d\,Im\,\zeta\over \zeta(\zeta-\lambda)},\ \
(\lambda,p)\in\Lambda^-_{\rho,\tau,\nu},&(4.15b)}$$
{\it where} $R_{\rho,\tau}$ {\it is defined by} (4.4).

\vskip 2 mm
{\bf Remark 4.1.}
In addition to (4.14), note that the definition of $(U_1,U_2)_{\rho,\tau}$
can be also written as
$$\eqalign{
&(U_1,U_2)_{\rho,\tau}(\lambda,p)=-{\pi\over 4}\int_{-\pi}^{\pi}
\biggl({|p|\over 2}
{{|\lambda|^2-1}\over \bar\lambda |\lambda|}(\cos\v-1)-{|p|\over \bar\lambda}
\sin\v\biggr)\times\cr
&U_1(z_1(\lambda,p,\v),-\xi(\lambda,p,\v))
U_2(z_2(\lambda,p,\v),p+\xi(\lambda,p,\v))\times\cr
&\chi_{2\tau\rho}(\xi(\lambda,p,\v))
\chi_{2\tau\rho}(p+\xi(\lambda,p,\v))d\v,\cr}\eqno(4.16)$$
where
$$\eqalign{
&z_1(\lambda,p,\v)={2k(\lambda,p)(\theta(-\xi(\lambda,p,\v))+i\omega
(-\xi(\lambda,p,\v)))\over i|p|},\cr
&z_2(\lambda,p,\v)={2(k(\lambda,p)+\xi(\lambda,p,\v))
(\theta(p+\xi(\lambda,p,\v))+i\omega
(p+\xi(\lambda,p,\v)))\over i|p|},\cr}\eqno(4.17)$$
$(\lambda,p)\in\Lambda_{\rho,\tau,\nu}$, $\v\in [-\pi,\pi]$,
$k(\lambda,p)$ is defined in (3.12b), $\xi(\lambda,p,\v)$ is defined by
(3.20), $\theta$, $\omega$ are the vector functions of (3.9), (3.10a).

We consider (4.11) as an integral equation for finding $H$ from $H^0$
with unknown remainder $Q_{\rho,\tau}$, where $H$, $H^0$, $Q_{\rho,\tau}$
are considered on $\Lambda_{\rho,\tau,\nu}$. Thus, actually, we consider
(4.11) as an approximate equation for finding $H$ on
$\Lambda_{\rho,\tau,\nu}$ from $H^0$ on $\Lambda_{\rho,\tau,\nu}$.
To deal with (4.11) we use Lemmas 4.1-4.5 given below.

Let
$$|||U|||_{\rho,\tau,\mu}=ess\,\sup\limits_{
(\lambda,p)\in\Lambda_{\rho,\tau,\nu}}(1+|p|)^{\mu}|U(\lambda,p)|
\eqno(4.18)$$
for $U\in L^{\infty}(\Lambda_{\rho,\tau,\nu})$, where $\rho>0$,
$\tau\in ]0,1[$, $\nu\in\S^2$, $\mu>0$.

\vskip 2 mm
{\bf Lemma 4.1.}
{\it Let} $v$ {\it and} $\rho$ {\it satisfy the same assumptions that in
Subsection  3.1. Let}  $R_{\rho,\tau}$, $Q_{\rho,\tau}$
 {\it be defined by} (4.4), (4.15), $\tau\in ]0,1[$. {\it Then}
$$\eqalignno{
&R_{\rho,\tau}\in L^{\infty}_{local}(\Lambda_{\rho,\tau,\nu}),&(4.19a)\cr
&|R_{\rho,\tau}(\lambda,p)|\le {3b(\mu_0,|\lambda|,|p|)C^2\over
(1-\eta)^2(1+2\tau\rho)^{\mu-\mu_0}(1+|p|)^{\mu_0}},\
(\lambda,p)\in\Lambda_{\rho,\tau,\nu},&(4.19b)\cr
&Q_{\rho,\tau}\in L^{\infty}(\Lambda_{\rho,\tau,\nu}),&(4.20a)\cr
&|||Q_{\rho,\tau}|||_{\rho,\tau,\mu_0}\le {3b_4(\mu_0)C^2\over
(1-\eta)^2(1+2\tau\rho)^{\mu-\mu_0}},&(4.20b)\cr}$$
{\it where} $2\le\mu_0\le\mu$, $b(\mu,|\lambda|,|p|)$ {\it is defined in}
(3.27), $\eta=\eta(C,\rho,\mu)$ {\it is defined by} (3.1),
$$b_4(\mu)={1\over \pi}(b_1(\mu)n_1+b_2(\mu)n_2+b_3(\mu)n_3),
\eqno(4.21)$$
{\it where} $b_1,b_2,b_3$ {\it are the constants of} (3.27) ({\it the
constants} $c_3,c_4,c_5$ {\it of Lemma} 5.2 {\it of} [No4]), $n_1,n_2,n_3$
{\it are the constants of Lemma} 11.1 {\it of} [No4].

Lemma 4.1 is proved in Section 7.

\vskip 2 mm
{\bf Lemma 4.2.}
{\it Let} $v$ {\it and} $\rho$ {\it satisfy the same assumptions that in
Subsection  3.1. Let}  $H^0$
 {\it be defined by} (4.8), (4.12), $\tau\in ]0,1[$. {\it Then}
$$\eqalignno{
&H^0\in L^{\infty}(\Lambda_{\rho,\tau,\nu}),&(4.22a)\cr
&|||H^0|||_{\rho,\tau,\mu_0}\le {C\over {1-\eta}}
\biggl(1+{b_4(\mu_0)C\over {1-\eta}}\biggr),&(4.22b)\cr}$$
{\it where} $2\le\mu_0\le\mu$, $\eta=\eta(C,\rho,\mu)$ {\it is defined by}
(3.1), $b_4$ {\it is defined by} (4.21).

Lemma 4.2 is proved in Section   7.

\vskip 2 mm
{\bf Lemma 4.3.}
{\it Let} $\rho>0$, $\nu\in\S^2$, $\tau\in ]0,1[$, $\mu\ge 2$. {\it Let}
$M_{\rho,\tau}$
 {\it be defined by} (4.13), (4.14) ({\it where} $\lambda,p$ {\it the
coordinates of Subsection} 3.3 {\it under assumption} (3.10a)). {\it Let}
$U_1,U_2\in L^{\infty}(\Lambda_{\rho,\tau,\nu})$,
$|||U_1|||_{\rho,\tau,\mu}< +\infty$, $|||U_2|||_{\rho,\tau,\mu}< +\infty$.
 {\it Then}
$$\eqalignno{
&M_{\rho,\tau}(U_j)\in L^{\infty}(\Lambda_{\rho,\tau,\nu}),\ j=1,2,
&(4.23)\cr
&|||M_{\rho,\tau}(U_j)|||_{\rho,\tau,\mu}\le c_8(\mu,\tau,\rho)
(|||U_j|||_{\rho,\tau,\mu})^2,\ j=1,2,&(4.24)\cr}$$
$$\eqalign{
&|||M_{\rho,\tau}(U_1)-M_{\rho,\tau}(U_2)|||_{\rho,\tau,\mu}\le\cr
&c_8(\mu,\tau,\rho)
(|||U_1|||_{\rho,\tau,\mu}+|||U_2|||_{\rho,\tau,\mu})
|||U_1-U_2|||_{\rho,\tau,\mu},\cr}\eqno(4.25)$$
{\it where}
$$c_8(\mu,\tau,\rho)=3b_1(\mu)\tau^2+4b_2(\mu)\rho^{-1}+4b_3(\mu)\tau,
\eqno(4.26)$$
{\it where} $b_1,b_2,b_3$ {\it are the constants of} (3.27).

Lemma 4.3 is proved in Section  7.

Lemmas 4.1, 4.2, 4.3 show that, under the assumptions of Proposition 4.1,
the nonlinear integral equation (4.11) for unknown $H$ can be analysed
for $H^0$, $Q_{\rho,\tau}$, $H\in L^{\infty}(\Lambda_{\rho,\tau,\nu})$
using the norm $|||\cdot|||_{\rho,\tau,\mu_0}$, where $2\le\mu_0\le\mu$.

Consider the equation
$$U=U^0+M_{\rho,\tau}(U),\ \ \rho>0,\ \ \tau\in ]0,1[, \eqno(4.27)$$
for unknown $U$ (where $U^0$, $U$ are functions on $\Lambda_{\rho,\tau,\nu}$).
Actually, under the assumptions of Proposition 4.1, we suppose that
$U^0=H^0+Q_{\rho,\tau}$ or consider $U^0$ as an approximation to
$H^0+Q_{\rho,\tau}$.

\vskip 2 mm
{\bf Lemma 4.4.}
{\it Let} $\rho>0$, $\nu\in\S^2$, $\tau\in ]0,1[$, $\mu\ge 2$
{\it and} $0<r< (2c_8(\mu,\tau,\rho))^{-1}$. {\it Let}
$M_{\rho,\tau}$
 {\it be defined by} (4.13), (4.14) ({\it where} $\lambda,p$ {\it are the
coordinates of Subsection} 3.3 {\it under assumption} (3.10a)). {\it Let}
$U^0\in L^{\infty}(\Lambda_{\rho,\tau,\nu})$ {\it and}
$|||U^0|||_{\rho,\tau,\mu}\le r/2$. {\it Then equation} (4.27) {\it is
uniquely solvable for}  $U\in L^{\infty}(\Lambda_{\rho,\tau,\nu})$,
$|||U|||_{\rho,\tau,\mu}\le r$, {\it and} $U$ {\it can be found by the
method of successive approximations, in addition,}
$$|||U-(M_{\rho,\tau,U^0})^n(0)|||_{\rho,\tau,\mu}\le
{r(2c_8(\mu,\tau,\rho)r)^n\over 2(1-2c_8(\mu,\tau,\rho)r)},\ n\in\N,
\eqno(4.28)$$
{\it where} $M_{\rho,\tau,U^0}$ {\it denotes the map}
$U\to U^0+M_{\rho,\tau}(U)$.

Lemma 4.4 is proved in Section 8  (using Lemma 4.3 and the lemma about
contraction maps).

\vskip 2 mm
{\bf Lemma 4.5.}
{\it Let the assumptions of Lemma} 4.4 {\it be fulfilled. Let also}
$\tilde U^0\in L^{\infty}(\Lambda_{\rho,\tau,\nu})$,
$|||\tilde U^0|||_{\rho,\tau,\mu}\le r/2$, {\it and} $\tilde U$ {\it
denote the solution of} (4.27) {\it with} $U^0$ {\it replaced by}
$\tilde U^0$, {\it where}
$\tilde U\in L^{\infty}(\Lambda_{\rho,\tau,\nu})$,
$|||\tilde U|||_{\rho,\tau,\mu}\le r$. {\it Then}
$$|||U-\tilde U|||_{\rho,\tau,\mu}\le
(1-2c_8(\mu,\tau,\rho)r)^{-1}|||U^0-\tilde U^0|||_{\rho,\tau,\mu}.
\eqno(4.29)$$

Lemma 4.5 is proved in Section 8.

Estimates (3.2), (3.3), Proposition 4.1 and Lemmas 4.1, 4.2, 4.3, 4.4, 4.5 
imply, in particular, the following result.

{\bf Proposition 4.2.}
{\it Let} $v$ {\it and} $\rho$ {\it satisfy the same assumptions that in
Subsection  3.1. Let}  $\nu\in\S^2$, $\tau\in ]0,1[$, $2\le\mu_0<\mu$.
{\it Let} $H$, $H^0$ {\it be defined on} $\Lambda_{\rho,\tau,\nu}$
{\it by} (4.8), (4.12) {\it and} $M_{\rho,\tau}$ {\it be defined by} (4.13).
{\it Let}
$$\eqalign{
&r_{min}(\mu,\mu_0,\tau,\rho,C)\le
r<(2c_8(\mu_0,\tau,\rho))^{-1},\cr
&r_{min}\buildrel \rm def \over = {2C\over {1-\eta(C,\rho,\mu)}}+
{2b_4(\mu_0)C^2\over (1-\eta(C,\rho,\mu))^2}
\biggl(1+{3\over (1+2\tau\rho)^{\mu-\mu_0}}\biggr),\cr}\eqno(4.30)$$
{\it where} $\eta$ {\it is defined in}
(3.1), $c_8$ {\it is defined by} (4.26). {\it Then the equation}
$$\tilde H_{\rho,\tau}=H^0+M_{\rho,\tau}(\tilde H_{\rho,\tau}) \eqno(4.31)$$
{\it is uniquely solvable for}
$\tilde H_{\rho,\tau}\in L^{\infty}(\Lambda_{\rho,\tau,\nu})$,
$|||\tilde H_{\rho,\tau}|||_{\rho,\tau,\mu_0}\le r$, {\it by the method of
successive approximations and}
$$|||H-\tilde H_{\rho,\tau}|||_{\rho,\tau,\mu_0}\le
{3b_4(\mu_0)C^2\over (1-2c_8(\mu_0,\tau,\rho)r)(1-\eta(C,\rho,\mu))^2
(1+2\tau\rho)^{\mu-\mu_0}}.\eqno(4.32)$$

Note that (4.30) can be fulfilled if and only if
$$r_{min}(\mu,\mu_0,\tau,\rho,C)<(2c_8(\mu_0,\tau,\rho))^{-1}.\eqno(4.33)$$

Using the definitions of $\eta$ of (3.1) and $c_8$ of (4.26) one can see that:
$$\eqalign{
&{\rm conditions}\ \ (3.1)\ \ {\rm and}\ \ (4.33) \ \ {\rm are\
fulfilled},\cr
&{\rm if}\ \ C\le c_9(\mu,\mu_0,\rho,\tau)\ \ {\rm for\ appropriate\
positive}\ c_9,\cr}\eqno(4.34)$$
where $2\le\mu_0<\mu$, $\ln\,\rho\ge 2$, $0<\tau<1$. Using (4.34) one can
see that Proposition 4.2 gives a method for approximate finding $H$ on
$\Omega_{\rho,\tau}^{\infty}$ from $H$ on  $b\Omega_{\rho,\tau}$  with
estimate
(4.32), at least, for sufficiently small potentials $v$ in the sense
$\|\hat v\|_{\mu}<C$, $C\le c_9(\mu,\mu_0,\rho,\tau)$.
However, the main point is that Proposition 4.2 also contains a global
result, see considerations given below in this section.

Due to (4.26) we have that
$$c_8(\mu_0,\tau,\rho)\le\ep\ \ {\rm if}\ \ 0<\tau\le\tau(\ep,\mu_0),\ \
\rho\ge\rho(\ep,\mu_0) \eqno(4.35)$$
for any arbitrary small $\ep>0$ and appropriate sufficiently small
$\tau(\ep,\mu_0)\in ]0,1[$ and sufficiently great $\rho(\ep,\mu_0)$.
Using (4.35) and the definition  of $\eta$ of (3.1) we obtain that:
$$\eqalign{
&{\rm conditions}\ \ (3.1), (4.33)
\ \ {\rm are\ fulfilled\ and}\cr
&0\le\eta(C,\rho,\mu)<\delta,\  \
0\le 2c_8(\mu_0,\tau,\rho)\,r_{min}(\mu,\mu_0,\tau,\rho,C)<\delta,\cr
&{\rm if}\ \ 0<\tau\le\tau_1(\mu,\mu_0,C,\delta),\ \
\rho\ge\rho_1(\mu,\mu_0,C,\delta),\cr}\eqno(4.36)$$
where $\tau_1$ and $\rho_1$ are appropriate constants such that
$\tau_1\in ]0,1[$ is sufficiently small and $\rho_1$ is
sufficiently great, $2\le\mu_0<\mu$, $0<\delta<1$.

As a corollary of Proposition 4.2 and property (4.36), we obtain the
following result.

\vskip 2 mm
{\bf Corollary 4.1.}
{\it Let} $v$ {\it satisfy} (2.1) {\it and} $\|\hat v\|_{\mu}\le C$. {\it Let}
$$0<\tau\le\tau_1(\mu,\mu_0,C,\delta),\ \ \rho\ge\rho_1(\mu,\mu_0,C,\delta),
\eqno(4.37)$$
{\it where} $2\le\mu_0\le\mu$, $0<\delta<1$.
 {\it Then} $H$ {\it on} $b\Omega_{\rho,\tau}$
{\it determines via} (4.12), (4.31) {\it the approximation}
$\tilde H_{\rho,\tau}$ {\it to} $H$ {\it on} $\Omega_{\rho,\tau}^{\infty}$
{\it with the error estimate} (4.32) ({\it where} $r$ {\it can be taken,
for example, as} $r=r_{min}$ {\it of}
 (4.30)) {\it and, in particular, with}
$$|||H-\tilde H_{\rho,\tau}|||_{\rho,\tau,\mu_0}=O(\rho^{-(\mu-\mu_0)})\ \
{\it as}\ \ \rho\to +\infty. \eqno(4.38)$$

The constant $C$ can be arbitrary great in Corollary 4.1 and, therefore, the
result of Corollary 4.1 is global.

\vskip 2 mm
{\bf 5. Approximate finding $\hat v$ on ${\cal B}_{2\tau\rho}$ from
$\tilde H_{\rho,\tau}$ on $\Omega_{\rho,\tau}^{\infty}$}

Consider, first, $H$ on $\Omega_{\rho,\tau,\nu}^{\infty}$ in the coordinates
$k$, $p$
 as $H$ on $\Lambda_{\rho,\tau,\nu}$ in the coordinates $\lambda$, $p$
according to (4.8). If $\hat v$ satisfies (2.1), then formulas (3.4),
(4.8), (3.12b), (3.13) imply that
$$\eqalign{
&H(\lambda,p)\to\hat v(p)\ \ {\rm as}\ \ \lambda\to 0,\cr
&H(\lambda,p)\to\hat v(p)\ \ {\rm as}\ \ \lambda\to\infty,\cr}\eqno(5.1)$$
where $p\in {\cal B}_{2\tau\rho}\b {\cal L}_{\nu}$, $\tau\in ]0,1[$.

Consider now $\tilde H_{\rho,\tau}$ defined in Proposition 4.2. Under the
assumptions of Proposition 4.2, the following formulas hold:
$$\eqalignno{
&\tilde H_{\rho,\tau}(\lambda,p)\to\hat v^+(p,\rho,\tau)\ \ {\rm as}\ \
\lambda\to 0,&(5.2a)\cr
&\tilde H_{\rho,\tau}(\lambda,p)\to\hat v^-(p,\rho,\tau)\ \ {\rm as}\ \
\lambda\to\infty,&(5.2b)\cr}$$
where
$$\eqalign{
&\hat v^+(p,\rho,\tau)={1\over 2\pi i}
\int\limits_{{\cal T}^+_{\rho/|p|}}H(\zeta,p)
{d\zeta\over \zeta}-\cr
&{1\over \pi}\int\!\!\!\int\limits_{{\cal D}^+_{\rho/|p|}}
(\tilde H_{\rho,\tau},\tilde H_{\rho,\tau})_{\rho,\tau}(\zeta,p)
{d\,Re\,\zeta d\,Im\,\zeta\over \zeta},\cr}\eqno(5.3a)$$
$$\eqalign{
&\hat v^-(p,\rho,\tau)={1\over 2\pi i}\int\limits_{{\cal T}^-_{\rho/|p|}}
H(\zeta,p)
{d\zeta\over \zeta}+\cr
&{1\over \pi}\int\!\!\!\int\limits_{{\cal D}^-_{\rho/|p|}}
(\tilde H_{\rho,\tau},\tilde H_{\rho,\tau})_{\rho,\tau}(\zeta,p)
{d\,Re\,\zeta d\,Im\,\zeta\over \zeta},\cr}\eqno(5.3b)$$
where $p\in {\cal B}_{2\tau\rho}\b {\cal L}_{\nu}$, $\tau\in ]0,1[$,
$(\tilde H_{\rho,\tau},\tilde H_{\rho,\tau})_{\rho,\tau}$ is defined
by means of (4.14), (4.16), (4.17).

Formulas (5.2), (5.3) follow from (4.31), where
$|||\tilde H_{\rho,\tau}|||_{\rho,\tau,\mu_0}<\infty$, $\mu_0\ge 2$, formulas
(4.12), (4.13) and estimate (3.27).

Formulas (5.1), (5.2), (4.18) imply that
$$\|\hat v-\hat v^{\pm}(\cdot,\rho,\tau)\|_{2\tau\rho,\mu_0}\le
|||H-\tilde H_{\rho,\tau}|||_{\rho,\tau,\mu_0},\eqno(5.4)$$
where
$$\|w\|_{r,\mu}=
ess\ \sup\limits_{p\in {\cal B}_r\b {\cal L}_{\nu}}(1+|p|)^{\mu}
|w(p)|,\ \ \mu>0,\ r>0,\eqno(5.5)$$
and $\rho,\tau,\mu_0$ are the same that in Proposition 4.2. Thus, under
the assumptions of Proposition 4.2 (or under the assumptions of
Corollary 4.1),  formulas (5.2), (5.3), (5.4), (4.32) imply that $\hat v$ on
${\cal B}_{2\tau\rho}$ can be approximately determined from
$\tilde H_{\rho,\tau}$ on $\Omega_{\rho,\tau}^{\infty}$ as
$\hat v_{\pm}(\cdot,\rho,\tau)$ of (5.2), (5.3) and
$$\eqalign{
&\|\hat v-\hat v^{\pm}(\cdot,\rho,\tau)\|_{2\tau\rho,\mu_0}\ \ {\rm is\
smaller\ or\ equal\ than\ the\ right-hand\ side\ of}\ \ (4.32)\cr
&{\rm and,\ in\ particular,\ is}\ \ O(\rho^{-(\mu-\mu_0)})\ \ {\rm as}\ \
\rho\to +\infty.\cr}\eqno(5.6)$$

\vskip 2 mm
{\bf 6. Reconstruction of $v$ from $\Phi$}

In this section we summarize our global 3D reconstruction
$$\Phi\buildrel 1\over \to h\big|_{b\Theta_{\rho,\tau}}\buildrel 2\over \to
\hat v\big|_{{\cal B}_{2\tau\rho}}\buildrel 3\over \to v \eqno(6.1)$$
developed in [No1], [No2] and in Sections 4, 5 of the present work. See
formulas (1.22), (1.3)-(1.12), (2.4) for notations used in (6.1). In (6.1)
the numbers $\rho>0$ and $\tau\in ]0,1[$ are parameters. The reconstruction
(6.1) for fixed $\rho$ and $\tau$ is approximate on the steps 2 and 3.
The steps 1, 2, 3 of (6.1) consist in the following:

(1) To find $h\big|_{b\Theta_{\rho,\tau}}$ from $\Phi$ we use formulas
and equations (1.23)-(1.25). In addition, if $v$ is sufficiently close to
some known non-zero background potential $v_0$, then instead of (1.23)-(1.25)
one can use their advanced version of [No2] for improving the reconstruction
stability.

(2) To find $\hat v$ on ${\cal B}_{2\tau\rho}$ from $h$ on
$b\Theta_{\rho,\tau}$ (approximately but stably and with minimal
approximation error) we proceed as follows
$$\eqalign{
&h\big|_{b\Theta_{\rho,\tau}}
{\ }_{\scriptstyle \longrightarrow \atop {\scriptstyle (1.4)}}\ \
H\big|_{b\Omega_{\rho,\tau}}\ \
{\ }_{\scriptstyle \longrightarrow \atop {\scriptstyle (3.12)}}\ \
H\big|_{b\Lambda_{\rho,\tau,\nu}}\cr
&{\ }_{\scriptstyle \longrightarrow \atop {\scriptstyle (4.12)}}\ \
H^0\big|_{\Lambda_{\rho,\tau,\nu}}\
{\ }_{\scriptstyle \longrightarrow \atop {\scriptstyle (4.31)}}\ \
\tilde H_{\rho,\tau}\ \ {\rm on}\ \
\Lambda_{\rho,\tau,\nu}\cr
&{\ }_{\scriptstyle \longrightarrow \atop {\scriptstyle (5.3)}}\ \
\hat v^{\pm}(\cdot,\tau,\rho)\ \ {\rm on}\ \
{\cal B}_{2\tau\rho},\cr}\eqno(6.2)$$
where $\hat v^{\pm}(\cdot,\tau,\rho)$ approximates $\hat v$ on
${\cal B}_{2\tau\rho}$. See also formulas (2.4), (2.10), (2.11), (3.15), 
(3.18) concerning the sets  $b\Omega_{\rho,\tau}$,
$b\Lambda_{\rho,\tau,\nu}$, $\Lambda_{\rho,\tau,\nu}$ mentioned in (6.2). In
(6.2) the substep from $H^0$ to $\tilde H_{\rho,\tau}$ consists in solving
the nonlinear integral equation (4.31), whereas all other substeps are
given by explicit formulas.

(3) Finally, from $\hat v^{\pm}(\cdot,\tau,\rho)$ of (6.2) we find
$v^{\pm}(\cdot,\tau,\rho)$ by formula (2.8), where $v^{\pm}(\cdot,\tau,\rho)$
approximates $v$ on $\R^3$.

One can see that on its steps 1 and 3 reconstruction (6.1) is reduced to
results of [No1], [No2] and to the inverse Fourier transform, whereas (6.2)
is developed in the present work. Some rigorous results concerning (2.8),
(6.2),  were already summarized as Theorem 2.1 and Corollary 2.1 of
Section 2.

In addition, under the assumptions of Theorem 2.1, a more detailed version
of (2.7) is given by
$$\eqalign{
&|\hat v(p)-\hat v^{\pm}(p,\tau,\rho)|\le {q(\mu,\mu_0,\tau,\rho,C) C^2\over
(1+2\tau\rho)^{\mu-\mu_0}},\cr
&q(\mu,\mu_0,\tau,\rho,C)={3b_4(\mu_0)\over
(1-2c_8(\mu_0,\tau,\rho) r_{min}(\mu,\mu_0,\tau,\rho,C))
(1-\eta(C,\rho,\mu))^2},\cr}\eqno(6.3)$$
where $\eta$, $b_4$, $c_8$, $r_{min}$ are defined in (3.1), (4.21), (4.26),
(4.30). Estimate (2.7) with
$$c_5(\mu,\mu_0,\tau,\delta)={3b_4(\mu_0)\over
(1-\delta)^3(2\tau)^{\mu-\mu_0}} \eqno(6.4)$$
follows from (6.3) and (4.36).

Note that
$h\big|_{b\Theta_{\rho,\tau}}$, $H\big|_{b\Omega_{\rho,\tau}}$ and
$H\big|_{b\Lambda_{\rho,\tau,\nu}}$
represent the same function in different coordinates. For stability analysis
of (6.2) it is convenient to fix this function as
$H\big|_{b\Lambda_{\rho,\tau,\nu}}$. Under the assumptions of Theorem 2.1,
this function has, in particular, the following properties (see (3.2), (3.3)
and the proof of (4.22b)):
$$\eqalignno{
&H\in {\cal C}(b\Lambda_{\rho,\tau,\nu}),&(6.5)\cr
&\|H\|_{\rho,\tau,\mu_0}\le {C\over {1-\eta(C,\rho,\mu)}},&(6.6)\cr
&\|T_bH\|_{\rho,\tau,\mu_0}\le {C\over {1-\eta(C,\rho,\mu)}}
\bigl(1+{b_4(\mu_0) C\over {1-\eta(C,\rho,\mu)}}\bigr),&(6.7)\cr}$$
where ${\cal C}$ denotes the space of continuous functions, $\eta$ is
defined by (3.1),
$$\eqalign{
&(T_bU)(\lambda,p)={1\over 2\pi i}\int\limits_{{\cal T}^+_{\rho/|p|}}
U(\zeta,p){d\zeta\over {\zeta-\lambda(1-0)}},\ \
(\lambda,p)\in b\Lambda^+_{\rho,\tau,\nu},\cr
&(T_bU)(\lambda,p)=-{1\over 2\pi i}\int\limits_{{\cal T}^-_{\rho/|p|}}
U(\zeta,p){\lambda d\zeta\over \zeta(\zeta-\lambda(1+0))},\ \
(\lambda,p)\in b\Lambda^-_{\rho,\tau,\nu},\cr}\eqno(6.8)$$
$$\|U\|_{\rho,\tau,\mu_0}=
\sup\limits_{(\lambda,p)\in b\Lambda_{\rho,\tau,\nu}}
(1+|p|)^{\mu_0}|U(\lambda,p)|,\eqno(6.9)$$
where $U$ is a test function on  $b\Lambda_{\rho,\tau,\nu}$,
$b\Lambda^{\pm}_{\rho,\tau,\nu}$ are defined in (4.9)

\noindent
(and $b\Lambda_{\rho,\tau,\nu}=
b\Lambda^+_{\rho,\tau,\nu}\cap b\Lambda^-_{\rho,\tau,\nu}$, see (4.10)).

Properties (6.5)-(6.7) are necessary properties of
$H\big|_{b\Lambda_{\rho,\tau,\nu}}$ under the assumptions of Theorem 2.1. In
addition, if two functions $H_1$, $H_2$ satisfy (6.5)-(6.7), where
$\tau$, $\rho$ satisfy (2.6), $C>0$, $2\le\mu_0<\mu$, $0<\delta<1$, then
$\hat v_i^{\pm}(\cdot,\tau,\rho)$ on ${\cal B}_{2\tau\rho}$ can be
constructed from $H_i$ via (4.12), (4.31), (5.2), (5.3), $i=1,2$ (in the same
way as $\hat v^{\pm}(\cdot,\tau,\rho)$ is constructed from
$H\big|_{b\Lambda_{\rho,\tau,\nu}}$ in the framework of Theorem 2.1), and
$$\|\hat v_1^{\pm}(\cdot,\tau,\rho)-
\hat v_2^{\pm}(\cdot,\tau,\rho)\|_{2\tau\rho,\mu_0}\le
(1-\delta)^{-1}\|T_b(H_1-H_2)\|_{\rho,\tau,\mu_0},\eqno(6.10)$$
where $\|\cdot\|_{2\tau\rho,\mu_0}$ in the left-hand side of (6.10) is
defined as in (5.5), $\|\cdot\|_{\rho,\tau,\mu_0}$ in the right-hand side of
(6.10) is defined by (6.9), $T_b$ is defined by (6.8).

The stability estimate (6.10) follows from:

(a) the maximum principle in $\lambda$ for $H_1^0$, $H_2^0$, $H_1^0-H_2^0$,
where $H_n^0$ is constructed from $H_n$ via (4.12), $n=1,2$,

(b) Lemma 4.5 and statement of (4.36),

(c) arguments similar with the arguments used for (5.4).

In the present work, restrictions in time prevent us from discussing the
stability of (6.1), (6.2) in more detail.

\vskip 2 mm
{\bf 7. Proofs of Lemmas 4.1, 4.2, 4.3}

\vskip 2 mm
{\it Proof of Lemma 4.1.}
Under the assumptions of Subsection 3.1, due to (3.2), (3.3), (3.25), (4.2),
we have that
$$\eqalign{
&H,\ \chi_{2\tau\rho}H,\ (1-\chi_{2\tau\rho})H\in
L^{\infty}_{\mu_0}(\Omega\b\bar\Omega_{\rho}),\cr
&|||H|||_{\rho,\mu_0}\le (1-\eta)^{-1}C,\
|||\chi_{2\tau\rho}H|||_{\rho,\mu_0}\le (1-\eta)^{-1}C,\cr
&|||(1-\chi_{2\tau\rho})H|||_{\rho,\mu_0}\le (1-\eta)^{-1}
(1+2\tau\rho)^{-(\mu-\mu_0)}C,\cr}\eqno(7.1)$$
where $\eta$ is given by (3.1), $\tau\in ]0,1[$, $0\le\mu_0\le\mu$.

Formulas (4.19) follow from (4.4), (7.1), (3.26), (3.27).

Formulas (4.20) follow from (4.15), (4.19), Lemma 11.1 of [No4] and the
following formulas
$$\eqalign{
&\int\limits_{{\cal D}^-_r}u_j(\zeta,s){|\lambda|\over |\zeta|}
{d\,Re\,\zeta\ d\,Im\,\zeta\over |\zeta-\lambda|}=
\int\limits_{{\cal D}^+_r}u_j(\zeta,s)
{d\,Re\,\zeta\ d\,Im\,\zeta\over |\zeta-\lambda^{-1}|},\cr
&j=1,2,3,\ \ r>1/2,\ \ \lambda\in {\cal D}^-_r,\ \ s>0,\cr}\eqno(7.2)$$
$$\eqalign{
&u_1(\zeta,s)={|\zeta|\over (|\zeta|^2+1)^2},\ \
u_2(\zeta,s)={(|\zeta|^2+1)s\over |\zeta|^2(1+s(|\zeta|+|\zeta|^{-1}))^2},\cr
&u_3(\zeta,s)={s\over |\zeta|(1+s(|\zeta|+|\zeta|^{-1}))}.\cr}\eqno(7.3)$$

Lemma 4.1 is proved.

\vskip 2 mm
{\it Proof of Lemma 4.2.}
Using (4.5), (4.12) and (4.8), (3.22) we obtain that
$$\eqalignno{
&H^0(\lambda,p)=H(\lambda,p)+{1\over \pi}
\int\!\!\!\int\limits_{{\cal D}^+_{\rho/|p|}}\{H,H\}(\zeta,p)
{d\,Re\,\zeta\ d\,Im\,\zeta\over {\zeta-\lambda}},\ \
(\lambda,p)\in\Lambda^+_{\rho,\tau,\nu},&(7.4a)\cr
&H^0(\lambda,p)=H(\lambda,p)+{1\over \pi}
\int\!\!\!\int\limits_{{\cal D}^-_{\rho/|p|}}\{H,H\}(\zeta,p)
{\lambda\ d\,Re\,\zeta\ d\,Im\,\zeta\over \zeta(\zeta-\lambda)},\ \
(\lambda,p)\in\Lambda^-_{\rho,\tau,\nu}.&(7.4b)\cr}$$

Formulas (4.22) follow from (7.4), (3.2), (3.3), (3.26), (3.27), Lemma 11.1 of
[No4] and formulas (7.2), (7.3).

Lemma 4.2 is proved.

\vskip 2 mm
{\it Proof of Lemma 4.3.}
Consider
$$\eqalignno{
&I_{\rho,\tau}(U,V)(\lambda,p)=I^+_{\rho,\tau}(U,V)(\lambda,p)=\cr
&-{1\over \pi}
\int\!\!\!\int\limits_{{\cal D}^+_{\rho/|p|}}(U,V)_{\rho,\tau}(\zeta,p)
{d\,Re\,\zeta\ d\,Im\,\zeta\over {\zeta-\lambda}},\ \
(\lambda,p)\in\Lambda^+_{\rho,\tau,\nu},&(7.5a)\cr
&I_{\rho,\tau}(U,V)(\lambda,p)=I^-_{\rho,\tau}(U,V)(\lambda,p)=\cr
&-{1\over \pi}
\int\!\!\!\int\limits_{{\cal D}^-_{\rho/|p|}}(U,V)_{\rho,\tau}(\zeta,p)
{\lambda\ d\,Re\,\zeta\ d\,Im\,\zeta\over \zeta(\zeta-\lambda)},\ \
(\lambda,p)\in\Lambda^-_{\rho,\tau,\nu},&(7.5b)\cr}$$
where $U,V\in L^{\infty}(\Lambda_{\rho,\tau,\nu})$,
$|||U|||_{\rho,\tau,\mu}< +\infty$,  $|||V|||_{\rho,\tau,\mu}< +\infty$,
$(U,V)_{\rho,\tau}$ is defined by (4.14), (4.16), $\rho$, $\nu$, $\tau$,
$\mu$ are the same that in Lemma 4.3.

Note that
$$\eqalignno{
&M_{\rho,\tau}(U_j)=I_{\rho,\tau}(U_j,U_j),\ \ j=1,2, &(7.6)\cr
&M_{\rho,\tau}(U_1)-M_{\rho,\tau}(U_2)=I_{\rho,\tau}(U_1-U_2,U_1)+
I_{\rho,\tau}(U_2,U_1-U_2),&(7.7)\cr}$$
where $U_1$, $U_2$ are the functions of Lemma 4.3.

Using (7.6), (7.7) one can see that in order to prove Lemma 4.3 it is
sufficient to prove that
$$\eqalignno{
&I_{\rho,\tau}(U,V)\in L^{\infty}(\Lambda_{\rho,\tau,\nu}),&(7.8)\cr
&|||I_{\rho,\tau}(U,V)|||_{\rho,\tau,\mu}\le c_8(\mu,\tau,\rho)
|||U|||_{\rho,\tau,\mu} |||V|||_{\rho,\tau,\mu}&(7.9)\cr}$$
under the same assumptions that in (7.5).

Formulas (7.8), (7.9) follow from (7.5), (4.14), (4.16), (3.26), (3.27),
(7.2), (7.3) and the following estimates:
$$
\int\limits_{{\cal D}^+_{\rho/|p|}}u_j(\zeta,|p|)
{d\,Re\,\zeta\ d\,Im\,\zeta\over |\zeta-\lambda|}\le\left\{\matrix{
(3/4)\pi (|p|/\rho)^2\ &\ \ {\rm for}\ \ j=1\hfill\cr
4\pi/\rho\ &\ \ {\rm for}\ \ j=2\hfill\cr
2\pi|p|/\rho\ &\ \ {\rm for}\ \ j=3,\hfill\cr}\right.\eqno(7.10)$$
where $u_1,u_2,u_3$ are defined by (7.3), $0<|p|< 2\tau\rho$, $0<\tau< 1$,
$\lambda\in {\cal D}^+_{\rho/|p|}$.

In turn, estimates (7.10) follow from the estimates
$$\zeta\in {\cal D}^+_r\Rightarrow |\zeta|\le (2r)^{-1},\ \ r\ge 1/2,
\eqno(7.11)$$
$$
\int\limits_{|\zeta|\le\ep}u_j(\zeta,s)
{d\,Re\,\zeta\ d\,Im\,\zeta\over |\zeta-\lambda|}\le\left\{\matrix{
3\pi\ep^2\ &\ \ {\rm for}\ \ j=1\hfill\cr
8\pi\ep/(\ep+s)\ &\ \ {\rm for}\ \ j=2\hfill\cr
4\pi\ep\ &\ \ {\rm for}\ \ j=3,\hfill\cr}\right.\eqno(7.12)$$
where $0\le\ep\le 1$, $s>0$, $|\lambda|\le\ep$. (To obtain (7.10) we use
(7.11), (7.12) for $r=\rho/|p|$, $\ep=(2r)^{-1}$, $s=|p|$.)

\vskip 2 mm
{\it The proof of (7.11).}
One can see that
$$\eqalign{
&\zeta\in  {\cal D}^+_r \buildrel (4.6) \over \Rightarrow
|\zeta|^2-4|\zeta|r+1>0,\ \ |\zeta|<1\Rightarrow\cr
&|\zeta|< 2r(1-\sqrt{1-1/(2r)^2})\le (2r)^{-1},\cr}\eqno(7.13)$$
where $r\ge 1/2$.

\vskip 2 mm
{\it The proof of (7.12).}
We have that
$$\eqalign{
&\int\limits_{|\zeta|\le\ep}u_j(\zeta,s)
{d\,Re\,\zeta\ d\,Im\,\zeta\over |\zeta-\lambda|}\le\cr
&\biggl(\int\limits_{|\zeta|\le\ep,\ |\zeta|\le |\zeta-\lambda|}+
\int\limits_{|\zeta|\le\ep,\ |\zeta|\ge |\zeta-\lambda|}\biggr)u_j(\zeta,s)
{d\,Re\,\zeta\ d\,Im\,\zeta\over |\zeta-\lambda|}\le A_j+B_j,\cr
&A_j=\int\limits_{|\zeta|\le\ep}u_j(\zeta,s)
{d\,Re\,\zeta\ d\,Im\,\zeta\over |\zeta|},\cr
&B_j=\int\limits_{|\zeta-\lambda|\le|\zeta|\le\ep}u_j(\zeta,s)
{d\,Re\,\zeta\ d\,Im\,\zeta\over |\zeta-\lambda|},\cr}\eqno(7.14)$$
where $j=1,2,3$. Further,
$$\eqalign{
&A_1\le\int\limits_{|\zeta|\le\ep}
d\,Re\,\zeta\ d\,Im\,\zeta=\pi\ep^2,\cr
&B_1\le\ep\int\limits_{|\zeta-\lambda|\le\ep}
{d\,Re\,\zeta\ d\,Im\,\zeta\over |\zeta-\lambda|}=2\pi\ep^2,\cr}\eqno(7.15)$$
$$\eqalign{
&A_2\le\int\limits_{|\zeta|\le\ep}
{(1+\ep^2)s\ d\,Re\,\zeta\ d\,Im\,\zeta
\over (|\zeta|+s(|\zeta|^2+1))^2|\zeta|},\cr
&B_2\le\int\limits_{|\zeta-\lambda|\le\ep}
{(1+\ep^2)s\ d\,Re\,\zeta\ d\,Im\,\zeta
\over (|\zeta-\lambda|+s(|\zeta-\lambda|^2+1))|\zeta-\lambda|},\cr
&A_2+B_2\le 4s\int\limits_{|\zeta|\le\ep}
{d\,Re\,\zeta\ d\,Im\,\zeta
\over (|\zeta|+s)^2|\zeta|}=8\pi s\int\limits_0^{\ep}{dr\over (r+s)^2}=
8\pi\ep/(\ep+s),\cr}\eqno(7.16)$$
$$\eqalign{
&A_3\le\int\limits_{|\zeta|\le\ep}
{s\ d\,Re\,\zeta\ d\,Im\,\zeta
\over (|\zeta|+s(|\zeta|^2+1))|\zeta|},\cr
&B_3\le\int\limits_{|\zeta-\lambda|\le\ep}
{s\ d\,Re\,\zeta\ d\,Im\,\zeta
\over (|\zeta-\lambda|+s(|\zeta-\lambda|^2+1))|\zeta-\lambda|},\cr
&A_3+B_3\le 2\int\limits_{|\zeta|\le\ep}
{s\ d\,Re\,\zeta\ d\,Im\,\zeta
\over (|\zeta|+s)|\zeta|}=4\pi s\int\limits_0^{\ep}{dr\over {r+s}}=
4\pi s\ln\bigl(1+{\ep\over s}\bigr)\le 4\pi\ep.\cr}\eqno(7.17)$$
Estimates (7.12) follow from (7.14)-(7.17).

Lemma 4.3 is proved.

\vskip 2 mm
{\bf 8. Proof of Lemmas 4.4 and 4.5}

\vskip 2 mm
{\it Proof of Lemma 4.4.}
For
$$\eqalign{
&U^0, U, U_1, U_2\in L^{\infty}(\Lambda_{\rho,\tau,\nu}),\cr
&|||U^0|||_{\rho,\tau,\mu}\le r/2,\ \  |||U|||_{\rho,\tau,\mu}\le r,\ \
|||U_1|||_{\rho,\tau,\mu}\le r,\ \ |||U_2|||_{\rho,\tau,\mu}\le r,\cr}
\eqno(8.1)$$
using Lemma 4.3 and the assumptions of Lemma 4.4 we obtain that
$$\eqalign{
&M_{\rho,\tau,U^0}(U)\in L^{\infty}(\Lambda_{\rho,\tau,\nu}),\cr
&|||M_{\rho,\tau,U^0}(U)|||_{\rho,\tau,\mu}\le
|||U^0|||_{\rho,\tau,\mu}+|||M_{\rho,\tau}(U)|||_{\rho,\tau,\mu}\le\cr
&r/2+c_8(\mu,\tau,\rho)r^2<r,\cr}\eqno(8.2)$$
$$\eqalign{
&|||M_{\rho,\tau,U^0}(U_1)-M_{\rho,\tau,U^0}(U_2)|||_{\rho,\tau,\mu}\le
2c_8(\mu,\tau,\rho) r |||U_1-U_2|||_{\rho,\tau,\mu},\cr
&2c_8(\mu,\tau,\rho) r<1,\cr}\eqno(8.3)$$
where
$$M_{\rho,\tau,U^0}(U)=U^0+M_{\rho,\tau}(U).\eqno(8.4)$$
Due to (8.1)-(8.4), $M_{\rho,\tau,U^0}$ is a contraction map of the ball
$U\in L^{\infty}(\Lambda_{\rho,\tau,\nu})$,
$|||U|||_{\rho,\tau,\mu}\le r$. Using now the lemma about contraction maps
and using the formulas
$$\eqalignno{
&|||U-M^n_{\rho,\tau,U^0}(0)|||_{\rho,\tau,\mu}\le
\sum_{j=n}^{\infty}|||M^{j+1}_{\rho,\tau,U^0}(0)-
M^j_{\rho,\tau,U^0}(0)|||_{\rho,\tau,\mu},&(8.5)\cr
&|||M_{\rho,\tau,U^0}(0)-0|||_{\rho,\tau,\mu}=|||U^0|||_{\rho,\tau,\mu}\le
r/2,&(8.6)\cr}$$
$$\eqalign{
&|||M^{j+1}_{\rho,\tau,U^0}(0)-M^j_{\rho,\tau,U^0}(0)|||_{\rho,\tau,\mu}
\buildrel (8.3)\over \le 2c_8(\rho,\tau,\mu) r\times\cr
&|||M^j_{\rho,\tau,U^0}(0)-M^{j-1}_{\rho,\tau,U^0}(0)|||_{\rho,\tau,\mu},\ \
j=1,2,3,\ldots,\cr}\eqno(8.7)$$
where $U$ is the  fixed point of $M_{\rho,\tau,U_0}$ in the aforementioned
ball, $M^0_{\rho,\tau,U^0}(0)=0$, we obtain Lemma 4.4.

\vskip 2 mm
{\it Proof of Lemma 4.5.}
We have that
$$\eqalignno{
&U-\tilde U=U^0-\tilde U^0+M_{\rho,\tau}(U)-M_{\rho,\tau}(\tilde U),&(8.8a)
\cr
&M_{\rho,\tau}(U)-M_{\rho,\tau}(\tilde U)\buildrel (7.7) \over =
I_{\rho,\tau}(U-\tilde U,U)+I_{\rho,\tau}(\tilde U,U-\tilde U),&(8.8b)\cr}$$
where  $I_{\rho,\tau}(U,V)$ is defined by (7.5).

In view of (8.8b) we can consider (8.8a) as a linear integral equation for
"unknown" $U-\tilde U$ with given $U^0-\tilde U^0$, $U$, $\tilde U$.
Using (7.9) and the properties $|||U|||_{\rho,\tau,\mu}\le r$,

\noindent
$|||\tilde U|||_{\rho,\tau,\mu}\le r$ we obtain that
$$|||I_{\rho,\tau}(U-\tilde U,U)-
I_{\rho,\tau}(\tilde U,U-\tilde U)|||_{\rho,\tau,\mu}\le 2c_8(\mu,\tau,\rho)
r |||U-\tilde U|||_{\rho,\tau,\mu}.\eqno(8.9)$$
Using (8.8b), (8.9) and solving (8.8a) by the method of successive
approximations we obtain (4.29). Lemma 4.5 is proved.

\vskip 4 mm
{\bf References}
\vskip 2 mm
\item{[  Al]} G.Alessandrini, {\it Stable determination of conductivity
by boundary measurements}, Appl. Anal. {\bf 27} (1988), 153-172.
\item{[ ABR]} N.V.Alexeenko, V.A.Burov and O.D.Rumyantseva, {\it Solution
of three-dimensional}

\item{      } {\it acoustical inverse scattering problem,II: modified
Novikov algorithm},  Acoust. J. {\bf 54}(3) (2008) (in Russian), English
transl.: Acoust. Phys. {\bf 54}(3) (2008).
\item{[  Am]} H.Ammari, {\it An Introduction to Mathematics of Emerging
Biomedical Imaging},

\item{      } Springer, Berlin, 2007.
\item{[  BC]} R.Beals and R.R.Coifman, {\it Multidimensional inverse
scattering and nonlinear partial differential equations}, Proc. Symp. Pure
Math. {\bf 43} (1985), 45-70.
\item{[ BRS]} V.A.Burov, O.D.Rumyantseva and T.V.Suchkova, {\it Practical
application possibilities of the functional approach to solving inverse
scattering problems}, Moscow: Physical Society, {\bf 3} (1990), 275-278
(in Russian).
\item{[   C]} A.-P.Calder\'on, {\it On an inverse boundary value problem},
Seminar on Numerical Analysis and its Applications to Continuum Physics
(Rio de Janeiro, 1980), pp.65-73, Soc. Brasil. Mat. Rio de Janeiro, 1980.
\item{[  F1]} L.D.Faddeev, {\it Growing solutions of the Schr\"odinger
equation}, Dokl. Akad. Nauk SSSR {\bf 165} (1965), 514-517 (in Russian);
English Transl.: Sov. Phys. Dokl. {\bf 10} (1966), 1033-1035.
\item{[  F2]} L.D.Faddeev, {\it Inverse problem of quantum scattering theory
II}, Itogi Nauki i Tekhniki, Sovr. Prob. Math. {\bf 3} (1974), 93-180
(in Russian); English Transl.: J.Sov. Math. {\bf 5} (1976), 334-396.
\item{[   G]} I.M.Gelfand, {\it Some problems of functional analysis and
algebra}, Proceedings of the International Congress of Mathematicians,
Amsterdam, 1954, pp.253-276.
\item{[  GN]} P.G.Grinevich, S.P.Novikov, {\it Two-dimensional "inverse
scattering problem" for negative energies and generalized-analytic functions.
I. Energies below the ground state}, Funkt. Anal. i Pril. {\bf 22}(1)
(1988), 23-33 (In Russian); English Transl.: Funkt. Anal. and Appl. {\bf 22}
(1988), 19-27.
\item{[  HM]} G.Henkin and V.Michel, {\it Inverse Conductivity Problem
on Riemann Surfaces},

\item{      } J.Geom.Anal. {\bf 18} (2008), 1033-1052.
\item{[  HN]} G.M.Henkin and R.G.Novikov, {\it The $\bar\pa$- equation in the
multidimensional inverse scattering problem}, Uspekhi Mat. Nauk {\bf 42(3)}
(1987), 93-152 (in Russian); English Transl.: Russ. Math. Surv. {\bf 42(3)}
(1987), 109-180.
\item{[   M]} N.Mandache, {\it Exponential instability in an inverse problem
for the Schr\"odinger equation}, Inverse Problems {\bf 17} (2001), 1435-1444.
\item{[ Na1]} A.I.Nachman, {\it Reconstructions from boundary measurements},
Ann. Math. {\bf 128} (1988), 531-576.
\item{[ Na2]} A.I.Nachman, {\it Global uniqueness for a two-dimensional
inverse boundary value problem}, Ann, Math. {\bf 142} (1995), 71-96.
\item{[ No1]} R.G.Novikov, {\it Multidimensional inverse spectral problem
for the equation $-\Delta\psi+(v(x)-Eu(x))\psi=0$}, Funkt. Anal. i Pril.
{\bf 22(4)} (1988), 11-22 (in Russian); English Transl.: Funct. Anal. and
Appl. {\bf 22} (1988), 263-272.
\item{[ No2]} R.G.Novikov, {\it Formulae and equations for finding scattering
data from the Dirichlet-to-Neumann map with nonzero background potential},
Inverse Problems {\bf 21} (2005), 257-270.
\item{[ No3]} R.G.Novikov, {\it The $\bar\pa$- approach to approximate inverse
scattering at fixed energy in three dimensions}, International Mathematics
Research Papers, {\bf 2005:6}, (2005), 287-349.
\item{[ No4]} R.G.Novikov, {\it On non-overdetermined inverse scattering at
zero energy in three dimensions}, Ann., Scuola Norm. Sup. Pisa Cl. Sci.
{\bf 5} (2006), 279-328.
\item{[  NN]} R.G.Novikov and N.N.Novikova, {\it On stable determination of
potential by boundary measurements}, ESAIM: Proceedings, to appear 2008.
\item{[  SU]} J.Sylvester and G.Uhlmann, {\it A global uniqueness theorem
for an inverse boundary value problem}, Ann. Math. {\bf 125} (1987),
153-169

\end